\documentclass[<options>]{article}
\usepackage[utf8]{inputenc}

\RequirePackage[numbers]{natbib}

\usepackage{graphicx,amsmath,amsfonts,amssymb,amsthm,psfrag,xcolor,framed, cancel, soul}

\usepackage{bbm,wrapfig,float}
\usepackage{fullpage}
\usepackage{comment}

\newtheorem{thm}{Theorem}[section]
\newtheorem{corollary}[thm]{Corollary}
\newtheorem{lemma}[thm]{Lemma}

\newtheorem{rem}[thm]{Remark}
\numberwithin{equation}{section}

\newtheorem{notation}[thm]{Notation}

\begin{document}

\title{Some exact results on Lindley process with Laplace jumps}
\author{Emanuele Lucrezia, Laura Sacerdote and Cristina Zucca\\
\small{University of Torino, Via Carlo Alberto 10,
10123 Turin, Italy,}\\ 
\small{emanuele.lucrezia@unito.it}\\
\small{laura.sacerdote@unito.it}\\
\small{cristina.zucca@unito.it}\\
}
\date{}
\maketitle
\begin{abstract}
We consider a Lindley process with Laplace distributed space increments. We obtain closed form recursive expressions for the density function of the position of the process and for its first exit time distribution from the domain $[0,h]$. We illustrate the results in terms of the parameters of the process. The work is completed by an open source version of the software. 
\end{abstract}
\par\medskip

 \par\medskip

\section{Introduction}

Random walk attracted mathematicians' attention both for their analytical properties and for their role in modeling instances, since the first time in which random walk were mentioned by Pearson \cite{Pearson}. Nowadays, random walks are universally recognized as the simplest example of stochastic process and are used as simplified version of more complex models. For modeling purposes it is often necessary to introduce suitable boundary conditions, constraining the evolution of the process on specific regions. Typical conditions are absorbing or reflecting boundaries and the involved random walk may be in one or more dimensions. The recurrence of the states of the free or bounded process in one or more dimensions was subject of past and recent studies and results are available for the free random walk \cite{Feller} or in presence of specific boundaries \cite{Cygan}.

The distribution of a simple random walk at time $n, n=1,2,\dots$ has a simple expression while, in the case of the distribution of the random walk  constrained by two absorbing boundaries, it involves important calculation efforts \cite{Cox}. 

Classical methods for the study of random walks refer to their Markov property. This happens, for example, when the focus is on recurrence properties in one or more dimensions. Alternatively, some results can be obtained by means of diffusion limits that give good approximations under specific hypotheses \cite{Lawler,Karlin,Kurtz81,Kurtz86} or allow to determine asymptotic results for specific problems such as the the maximum of a random walk or the first exit time across general boundaries \cite{Gut}. However, the presence of boundaries often determines combinatorics problems and discourages the search of closed form solutions. Furthermore, the switch from unitary jumps to continuous jumps introduces important computational difficulties. There are thousands of papers about random walk and it is impossible to refer to all of them. We limit ourselves to cite the recent excellent review by Dshalalow and White \cite{Dshalalow} that lists many of the most important available results.

In this paper we focus on a particular constrained random walk: the Lindley process \cite{Lindley}. In particular, we prove a set of analytical results about it. This process is a discrete time random walk characterized by continuous jumps with a specified distribution. It was introduced in \cite{Lindley} to describe waiting times experienced by customers in a queue over time.
In particular, it is defined as 
\begin{eqnarray} \label{Lindley0}
W_n&=&\max(0,W_{n-1}+Z_{n})\\
W_0&=&x\geq 0.\nonumber
\end{eqnarray}
where
\begin{itemize}
\item $T_n$ is the time between the $n$-th and $(n+1)$-th arrivals,
\item $S_n$ is the service time of the $n$-th customer,
\item $Z_n = S_n-T_n$
\item $W_n$ is the waiting time of the $n$-th customer.
\end{itemize}
 Modeling interest for the Lindley process is not limited to queuing problems. Indeed, this process arises in reliability context, in sequential test framework, through the study of CUSUM tests  \cite{Page}. Furthermore, the same process appears in problems related to resources management \cite{Bhattacharya1} and in highlighting atypical regions in biological sequences, transferring biological sequence analysis tools to break-point detection for on-line monitoring \cite{Mercier}.

 Many contributions concerning properties of the Lindley process are motivated by their important role in applications \cite{ Bhattacharya1, Bhattacharya, Dshalalow1, Feller, Lindley}.
 A large part of the papers about the Lindley process concerns its asymptotic behaviour. Using the strong law of large numbers, in 1952 Lindley \cite{Lindley}  showed that the process admits a limit distribution as $n$ diverges if and only if $\mathbb{E}[Z]<0$, i.e. if the customer's arrival rate is slower than the service one. Furthermore, when $\mathbb{E}[Z]=0$ the ratio $W_n/\sqrt{n}$ converges to the modulus of a Gaussian random variable.
 Lindley also showed that the limit distribution is solution of an integral equation and Kendall solved it in the case of exponential arrivals \cite{Kendall51} while Erlang determined its expression when both arrival and services times are exponential \cite{Brockmeyer}. Simulations were used by \cite{Iams} to determine invariant distribution when $Z_i, i=1, 2,  \cdots$ are Laplace or Gaussian distributed. The analytical expression of the limit distribution was determined by Stadje \cite{Stadje} for the case of integer valued increments. Recent contributions consider recurrence properties in higher dimensions \cite{Cygan, Diaconis, Peigne} or the study of the rate of convergence of expected values of functions of the process.
 
 Other contributions make use of stochastic ordering notion \cite{Lakatos} or of continuous time versions of the Lindley process \cite{Asmussen}. Recently, the use of machine learning techniques has been proposed to learn Lindley recursion (\ref{Lindley0}) directly from the waiting time data of a G/G/1 queue \cite{Palomo}. Furthermore, Lakatos and Zbaganu \cite{Lakatos}, and Raducan, Lakatos and Zbaganu \cite{Raducan} introduced a class of computable Lindley processes for which the distribution of the space increments is a mixture of two distributions. To the best of our knowledge no other analytical expressions are available for a Lindley process.

 Concerning the first exit time problem for the Lindley process, it has been considered mainly in the framework of CUSUM methodology to detect the time in which a change of the parameters of the process happens. In this context the focus was on the Run Length that corresponds to the exit time of the process from a boundary and in \cite{Gan, Wardeman} such distribution and its expected value is determined when the $Z_i, i=1,2, \dots$ are exponentially distributed with shift. Markovich and Razumchick \cite{Markovich} investigated a problem related to first exit times, i.e. the appearance of clusters of extremes defined as subsequent
expediences of high thresholds in a Lindley process.

Here we consider a Lindley process characterized by Laplace distributed space increments.  The use of such distribution for the increments is easily interpreted in the case of a queueing model while its interest in the CUSUM framework will be exemplified in Section \ref{Sect Lindley}. 
Our aim is to derive expressions for the distribution of the process as the time evolves together with its first exit time (FET) from $[0,h]$. we prove recursive closed form formulae for such distributions and we show that different formulae hold on different parameters domain. After the preliminary introduction of the process presented in Section 2, in Section 3 and 4 we study the distribution of the process and its first exit times from $[0,h]$, respectively. In these sections we do not present any proofs that we postpone in Section 5 and 6. Due to the complexity of the derived exact formulae for the studied distributions, we completed our work implementing the software necessary for a fast computation of the formulae of interest. This software is open source and can be found in the GitHub repository \cite{Lucrezia}.  Lastly, in Section 7 we illustrate the role of the parameters of the Laplace distribution on the position and FET distribution of the Lindley process.

\section{Lindley process with Laplace space increments}\label{Sect Lindley}
Let us consider a random walk with i.i.d. space increments $Z$ distributed as a Laplace random variable characterized by parameter $\mu$ and scale parameter $\sigma$, i.e. the density function of the random variable $Z$ is 
\begin{equation}\label{Laplace}
f_Z(z)= \frac{e^{-\frac{|z-\mu|}{\sigma}}}{2\sigma}
\end{equation}
and the mean and the variance are $\mathbb{E}(Z)=\mu$ and $V(Z)=2\sigma^2$.
The process we are interested in is $W=\{W_n,n\geq 0\}$ where
\begin{eqnarray} \label{Lindley}
W_n&=&\max(0,W_{n-1}+Z_{n})\\
W_0&=&x>0.\nonumber
\end{eqnarray}

\begin{figure}
\centering
\includegraphics[height= 8cm, keepaspectratio]{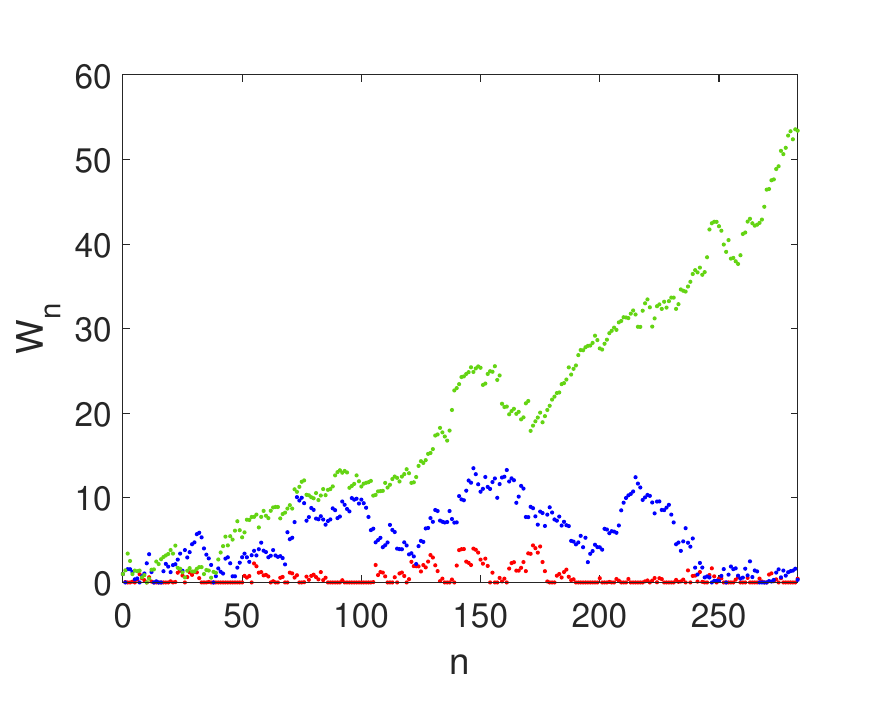}
\caption{Trajectories of the Lindley process with $\sigma=1$, $x=1$ and $\mu=-0.3$ (red), $\mu=0$ (blue) and $\mu=0.3$ (green).}\label{Fig posizione 1}
\end{figure}
In Figure \ref{Fig posizione 1} are plotted three trajectories of the Lindley process with $\sigma=1$, $x=1$ and different $\mu$.

As in the case of the simple random walk, here the behaviour of the random walk dramatically changes according to the sign of $\mu$. Indeed, it is well known \cite{Gut} that, when $\mu>0$, all the states are transient and the process has a drift toward infinity. While, in the case $\mu=0$ the process is null recurrent  and,  when $\mu<0$, the process is positive recurrent and admits limit distribution.
Such distribution is solution of the Lindley integral equation \cite{Lindley}. Unfortunately the analytical solution is unknown and numerical methods could be applied. 

The recurrence properties imply certainty of the first exit time of the process $W$ from $[0,h]$, for each $h>x$, and the study of the first exit time distribution can be performed.

The distribution of the process $W$ is not known and will be the subject of Chapter 3  while the distribution of the first exit time will be the subject of Chapter 4.  

An example in which the Lindley process with Laplace increments plays a role in the CUSUM framework is presented here. Suppose we have a sequence $\{X_n\}_{n \ge 1}$ of independent random variables with probability density function that changes at an unknown time $m$. The pre-change density is $f(x)$ and the post-change one is $g(x)$, that is 
\begin{equation*}
	X_n \sim  \left\{
	\begin{array}{lr}
		f(x),\hspace{0.3cm} \text{if } 1\le n < m \\
		g(x),\hspace{0.3cm} \text{if } n \ge m 
	\end{array}
	\right.
\end{equation*}
where $f(x)$ is a Laplace density function (\ref{Laplace}) 
and 
\begin{equation}
    g(x)=f(x)e^{\theta x-b(\theta)}
\end{equation} 
with $\theta$ positive parameter and \begin{equation}
    b(\theta)=\log(\mathbbm{E}[e^{\theta X}])=\mu \theta -\log(1-\sigma^2\theta^2)
\end{equation} with $X\sim f(x)$ and $-1/\sigma<\theta< 1/\sigma$.
This means that, up to time $m$, the random variable $X_n$ is Laplace distributed $X_n\sim Laplace(\mu, \sigma)$ with mean $\mu$ while, after time $m$ it corresponds to a skewed Laplace distribution with mean  $\mu+(2\sigma^2\theta)/(1-\sigma^2\theta^2)>\mu$. 

In this case that the loglikelihood ratio of the $n$-th observation is  
\begin{equation}
    LLR_n=\log\left(\frac{g(X_n)}{f(X_n)}\right)=\theta X_n-b(\theta),
\end{equation} i.e. it is a linear function of $X_n$ with slope $\theta$ and specific intercept for each slope.  The distribution of the loglikelihood ratio is then a rescaled and shifted Laplace random variable  $X_n\sim Laplace(\theta \mu-b(\theta), \sigma \theta )$. The knowledge of the distribution of the position and of the exit times of such distribution presented in the following sections allows to perform the test avoiding approximations.

\section{Distribution of $W_n$.}\label{Sect. distribution}

Let us consider the Lindley process (\ref{Lindley}) with Laplace distributed jumps (\ref{Laplace}).
In the following we will use the distribution of the process at time $n \geq 0$ that we denote as
\begin{equation}
F_n(u|x):=F_{W_n}(u|x)=\mathbb{P}(W_n \leq u|W_0=x).
\end{equation}
and we indicate with 
\begin{equation}
f_n(u|x):=f_{W_n}(u|x)=\frac{\partial}{\partial u}\mathbb{P}(W_n \leq u|W_0=x).
\end{equation}
the corresponding probability density function, where the derivative is the distributional derivative, because the $W_n$ are mixed random variables. 

In the following we make use of Dirac delta function $\delta(u)$ with the agreement that, for each $x>0$ 
\begin{equation*}
    \int_0^x f(u)\delta(u)=f(0)
\end{equation*}
Since $W_n \geq 0$, the density $f_n(u|x)=0$ for $u< 0$. 
Moreover, if $a>b$ the sum $\sum_{i=a}^b f_i=0$.

\begin{lemma} \label{Lemma1}
The probability distribution function of $W_1$ for $\mu>-x$ is
\begin{align}\label{F1(u|x)}
F_1(u|x)=F_{W_1}(u|x)=\left\lbrace 
\begin{array}{ll}
0 & u<0\\
\frac{e^{-\frac{\mu+x}{\sigma}}}{2}  &u=0\\
\frac{1}{2}e^{-\frac{\mu-u+x}{\sigma}}  &0<u \leq x+\mu\\
1-\frac{1}{2}e^{\frac{\mu-u+x}{\sigma}} & u>x+\mu
\end{array}
\right.
\end{align}
while, for  $\mu\leq -x$ is
\begin{align}\label{F1(u|x) mu<-x}
F_1(u|x)=F_{W_1}(u|x)=\left\lbrace 
\begin{array}{ll}
0 & u<0\\
1-\frac{1}{2}e^{\frac{\mu+x}{\sigma}}  &u=0\\
1-\frac{1}{2}e^{\frac{\mu-u+x}{\sigma}} & u>0
\end{array}
\right.
\end{align}   
The corresponding probability density function is 
\begin{align}\label{f_1 generale}
 f_1(u|x)=f_{W_1}(u|x)=\left\{
 \begin{array}{ll}
       \frac{e^{-\frac{|u-x-\mu|}{\sigma}}}{2\sigma} \mathbbm{1}_{(0,+\infty)}(u)+\frac{e^{-\frac{(x+\mu)}{\sigma}}}{2} \delta(u)     & x>-\mu\\
       \frac{e^{-\frac{u-x-\mu}{\sigma}}}{2\sigma} \mathbbm{1}_{(0,+\infty)}(u)+\left(1-\frac{e^{\frac{(x+\mu)}{\sigma}}}{2}\right) \delta(u)   &  x\leq -\mu
 \end{array}
  \right.
 \end{align}
where $\delta(u)$ is the Dirac delta function. 
\end{lemma}
\begin{rem}
Using the Markov property of the process $W$, the one-step transition probability density function is 
\begin{align} \label{denstransiz}
    f(u,n+1|y,n)&:=\frac{\partial}{\partial u}\mathbb{P}(W_{n+1}\leq u|W_n=y)=f_1(u|y)\\
    &=\left\{
 \begin{array}{ll}
       \frac{e^{-\frac{|u-y-\mu|}{\sigma}}}{2\sigma} \mathbbm{1}_{(0,+\infty)}(u)+\frac{e^{-\frac{(y+\mu)}{\sigma}}}{2} \delta(u)     & y>-\mu\\
       \frac{e^{-\frac{u-y-\mu}{\sigma}}}{2\sigma} \mathbbm{1}_{(0,+\infty)}(u)+\left(1-\frac{e^{\frac{(y+\mu)}{\sigma}}}{2}\right) \delta(u)   &  y\leq -\mu
 \end{array}
  \right. \nonumber
\end{align}
\end{rem}
\begin{notation}
In the forthcoming, when not necessary we will skip the dependence on the initial position of the process: $F_n(u|x)=F_n(u)$, $f_n(u|x)=f_n(u)$.
\end{notation}

In order to determine the distribution of $W_n$ computations change according to the sign of $\mu$. Theorem \ref{Theorem Position mu>=0} gives the distribution for $\mu>0$. When $\mu$ is negative two different cases arise when $-x<\mu$ (Theorem \ref{Theorem Position -x<mu<0}) and $\mu>-x$ (Theorem \ref{Theorem Position mu<=-x}).

\begin{thm} \label{Theorem Position mu>=0}
For a Lindley process $\{W_n, n\geq 0\}$ characterized by Laplace increments with location parameter $\mu\geq0$, the probability density function of the position is given by 
\begin{subequations}\label{fn}
\begin{align}
f_{n}(u)&=\sum _{i=0}^{n+1} f_n^i(u), \qquad \qquad u\geq 0, n\geq1\\
f_{0}(u)&=\delta (u-x)
\end{align}
\end{subequations} 
where 
\begin{subequations}
\begin{align}
&f_n^i(u)=\left\{\sum_{j=0}^{m^i_n-1} \frac{1}{(2 \sigma)^n}\left[a_{n}^{(i,j)}(u-n\mu)^j e^{\frac{u}{\sigma}}+b_{n}^{(i,j)}(u-n\mu)^j e^{-\frac{u}{\sigma}}\right]\right\} \mathbbm{1}_{I_n^i}(u) \qquad 1 \leq i\leq n+1 \label{fniPos}\\
&f_n^0(u)=c_n \delta(u)=c_n  \mathbbm{1}_{I_n^0}(u)\label{fdirac}
\end{align}
\end{subequations}
$ m_n^i=\min(n,i)$ and $I_n^i$, $i=0, \dots n+1$, with $\bigcup_{i=0}^{n+1} I_n^i=[0, +\infty)$, are
\begin{subequations}\label{I_n}
\begin{align}
&I_n^0=\{0\}\\
&I^i_n=((i-1)\mu,i\mu] \qquad i=1,\dots , n-1\\
&I_n^n=((n-1)\mu,n\mu+x]\\
&I_n^{n+1}=(n \mu+x, +\infty),
\end{align}
\end{subequations}
The coefficients $a_{n+1}^{(i,j)}$, $b_{n+1}^{(i,j)}$ and $c_{n+1}$ verify the following recursive relations for $1\leq i \leq n+1$, $1\leq j \leq m_n^i-1$
\begin{subequations} \label{recursionPos}
\begin{align}
a_{n+1}^{(i,j)}&=e^{-\frac{\mu}{\sigma}}\left[\sum_{k=j-1}^{m^{i-1}_n-1} (-1)^{k+j} a_{n}^{((i-1),k)}
\frac{k!}{j!}   \left(\frac{2}{\sigma}\right)^{k-j+1}
\right]\\
b_{n+1}^{(i,j)}&=-e^{\frac{\mu}{\sigma}}\left[\sum_{k=j-1}^{m^{i-1}_n-1} b_{n}^{((i-1),k)}\frac{k!}{j!} \left(\frac{2}{ \sigma}\right)^{k-j+1}\right]  \\
a_{n+1}^{(i,0)}&=e^{-\frac{\mu}{\sigma}}\left[\sum _{j=i}^{n+1}  A^{(2)}_j +\sum_{k=0}^{m^{i-1}_n-1}\left[a_{n}^{((i-1),k)}\left(  \frac{((i-1-n)\mu)^{k+1}}{k+1}+(-1)^{k}k!\left(\frac{\sigma}{2}\right)^{k+1} \right)\right.\right.\nonumber\\
&\left.-b_{n}^{((i-1),k)}\frac{e^{\frac{-2n\mu}{\sigma}}\sigma^{k+1}}{2^{k+1}}\left[\Gamma\left(1+k,\frac{2((i-1-n)\mu)}{\sigma}\right)\right]+c_n\delta_{i,1}\right]\\
b_{n+1}^{(i,0)}&=e^{\frac{\mu}{\sigma}}\left[\sum _{j=1}^{i-2}  A^{(1)}_j- \sum_{k=0}^{m^{i-1}_n-1}
\left[b_{n}^{((i-1),k)} \left( \frac{((i-2-n)\mu)^{k+1}}{k+1}-k!\left(\frac{\sigma}{2}\right)^{k+1} \right)
\right]\right. \nonumber\\
& \left.-a_{n}^{((i-1),k)}\frac{e^{\frac{2n\mu}{\sigma}}\sigma^{k+1} (-1)^k}{2^{k+1}}\left[\Gamma\left(1+k,\frac{-2((i-1-n)\mu)}{\sigma}\right) \right]+c_n(1-\delta_{i,1})\right]\\
c_{n+1}&=\frac{e^{-\frac{\mu}{\sigma}}}{2(2 \sigma)^n}\left\{\sum _{j=1}^{n+1} \sum_{k=0}^{m^i_n-1}  \left[\left.a_{n}^{(j,k)}\frac{(y-n\mu)^{k+1}}{k+1}\right|_{I_n^j} -b_{n}^{(j,k)}\frac{e^{-\frac{2n\mu}{\sigma}}\sigma^{k+1}}{2^{k+1}}\left.\Gamma\left(1+k,\frac{2(y-n\mu)}{\sigma}\right)\right|_{I_n^j}\right]+(2 \sigma)^n c_n\right\} \label{c_n}
\end{align}
\end{subequations}
where $\delta_{i,j}$ is the Kronecker delta and 
\begin{subequations}\label{A_j}
\begin{align}
A^{(1)}_j&:=\left.A^{(1)}_j(y)\right|_{I_{n}^j}=\sum_{k=0}^{m^j_n-1}   \left\{a_{n}^{(j,k)}\frac{e^{\frac{2n\mu}{\sigma}}\sigma^{k+1} (-1)^k}{2^{k+1}}\left[\left.\Gamma\left(1+k,\frac{-2(y-n\mu)}{\sigma}\right)\right|_{I_{n}^j} \right]+b_{n}^{(j,k)}\left.\frac{(y-n\mu)^{k+1}}{k+1}\right|_{I_{n}^j} \right\} \\
A^{(2)}_j&:=\left.A^{(2)}_j(y)\right|_{I_{n}^j}=\sum_{k=0}^{m^j_n-1} \left\{a_{n}^{(j,k)}\left.\frac{(y-n\mu)^{k+1}}{k+1}\right|_{I_{n}^j} -b_{n}^{(j,k)} \frac{e^{\frac{-2n\mu}{\sigma}}\sigma^{k+1}}{2^{k+1}}\left[\left.\Gamma\left(1+k,\frac{2(y-n\mu)}{\sigma}\right)\right|_{I_{n}^j}\right]
\right\}.\label{A_j_2}
\end{align} 
\end{subequations}

The initial values for the recursion of (\ref{recursionPos}) are
\begin{subequations}\label{CI mu>0}
\begin{align}
&c_1=\frac{e^{-\frac{\mu+x}{\sigma}}}{2}\\
&a_{1}^{(1,0)}=e^{-\frac{\mu+x}{\sigma}}  \\
&b_{1}^{(1,0)}=a_{1}^{(2,0)}=0\\
&b_{1}^{(2,0)}=e^{\frac{\mu+x}{\sigma}}
\end{align} 
\end{subequations}

\end{thm}

\begin{rem}
Observe that the coefficients $a_n^{(n+1,k)}\equiv 0$, for each admissible $k$. This prevents the terms in (\ref{c_n}) and (\ref{A_j_2}) to explode.
\end{rem}

The following corollary may be useful for computational purposes.

    \begin{corollary}\label{Cor c_n+1}
   The constant coefficients $c_n$, $n=1,2,\dots$ can also be obtained as
 \begin{equation}\label{c f}
     c_{n+1}=\sigma f_{n+1}^1(0).
 \end{equation}
\end{corollary}

\begin{rem}
Observe that the density (\ref{fn}) refers to a mixed random variable.
\end{rem}

\begin{thm} \label{Theorem Position -x<mu<0}
For a  Lindley process $\{W_n, n\geq 0\}$ characterized by Laplace increments with location parameter $-x<\mu<0$, the probability density function of the position is given by 
 \begin{subequations}\label{fn_b}
 \begin{align}
f_{n}(u)&=\sum _{i=0}^{2} f_n^i(u), \qquad \qquad u\geq 0, n\geq 1\\
f_{0}(u)&=\delta (u-x)
\end{align}
\end{subequations}  
where 
 \begin{subequations}\label{f_n T2} 
\begin{align}
&f_n^i(u)=\left\{\sum_{j=0}^{n-1} \frac{1}{(2 \sigma)^n}\left[a_{n}^{(i,j)}(u-n\mu)^j e^{\frac{u}{\sigma}}+b_{n}^{(i,j)}(u-n\mu)^j e^{-\frac{u}{\sigma}}\right]\right\}\mathbbm{1}_{I_n^i}(u) \qquad i=1,2 \label{fni 3.4}\\
&f_n^0(u)=c_n \delta(u)= c_n\mathbbm{1}_{I_n^0}(u)\label{fdirac_a}
\end{align}
\end{subequations}  
and $I_n^i$, $i=0,1,2$ are
 \begin{subequations}\label{partitionTh2}
 \begin{align}
&I_n^0=\{0\}\\
&I_n^1=(0, \max(0,n\mu+x))\\
&I_n^{2}=(\max (0,n \mu+x), +\infty)
\end{align}
\end{subequations}
\begin{itemize}
    \item
\underline{If $x+n\mu>0$}, the coefficients $a_{n+1}^{(i,j)}$, $b_{n+1}^{(i,j)}$  verify the recursive relations for $i=1,2$ and $1\leq j\leq n$
 \begin{subequations}\label{coeff T2}
\begin{align}
a_{n+1}^{(i,j)}&=e^{-\frac{\mu}{\sigma}}\left[\sum_{k=j-1}^{n-1} (-1)^{k+j} a_{n}^{(i,k)}
\frac{k!}{j!}   \left(\frac{\sigma}{2}\right)^{k-j+1}
\right]\label{ricorsione a}\\
b_{n+1}^{(i,j)}&=e^{\frac{\mu}{\sigma}}\left[\sum_{k=j-1}^{n-1} b_{n}^{(i,k)}\frac{k!}{j!} \left(\frac{\sigma}{2}\right)^{k-j+1}\right]\\
a_{n+1}^{(1,0)}&=e^{-\frac{\mu}{\sigma}}\left[ B_2 +\sum_{k=0}^{n-1}\left[a_{n}^{(1,k)}\left(  \frac{(x)^{k+1}}{k+1}+(-1)^{k}k!\left(\frac{\sigma}{2}\right)^{k+1} \right)\right.\right. \nonumber\\
&\left.-b_{n}^{(1,k)}\frac{e^{\frac{-2n\mu}{\sigma}}\sigma^{k+1}}{2^{k+1}}\left[\Gamma\left(1+k,\frac{2x}{\sigma}\right)\right]\right]\\
b_{n+1}^{(1,0)}&=e^{\frac{\mu}{\sigma}}\left[\sum_{k=0}^{n-1}
\left[b_{n}^{(1,k)} \left( -\frac{(-n\mu)^{k+1}}{k+1}+k!\left(\frac{\sigma}{2}\right)^{k+1} \right)
\right]\right.\nonumber\\
& \left.-a_{n}^{(1,k)}\frac{e^{\frac{2n\mu}{\sigma}}\sigma^{k+1} (-1)^k}{2^{k+1}}\left[\Gamma\left(1+k,\frac{2n\mu}{\sigma}\right) \right]+(2\sigma)^nf_n^0(0)\right]\\
a_{n+1}^{(2,0)}&=0
 \end{align}
 \begin{align}
b_{n+1}^{(2,0)}&=e^{\frac{\mu}{\sigma}}\left[B_1+\sum_{k=0}^{n-1}
\left[b_{n}^{(2,k)} \left( -\frac{x^{k+1}}{k+1}+k!\left(\frac{\sigma}{2}\right)^{k+1} \right)
\right]\right. \nonumber\\
& \left.-a_{n}^{(2,k)}\frac{e^{\frac{2n\mu}{\sigma}}\sigma^{k+1} (-1)^k}{2^{k+1}}\left[\Gamma\left(1+k,-\frac{2x}{\sigma}\right) \right]+(2\sigma)^nf_n^0(0)\right]
\end{align}
\end{subequations}
where 
\begin{subequations}\label{B_i}
 \begin{align} 
B_1&= \sum_{k=0}^{n-1} \left\{a_{n}^{(1,k)}\frac{e^{\frac{2n\mu}{\sigma}}\sigma^{k+1} (-1)^k}{2^{k+1}}\left[\left.\Gamma\left(1+k,\frac{-2(y-n\mu)}{\sigma}\right)\right|_{0}^{x+n\mu} \right]+b_{n}^{(1,k)}\left.\frac{(y-n\mu)^{k+1}}{k+1}\right|_{0}^{x+n\mu} \right\}\\
B_2&= \sum_{k=0}^{n-1} \left\{b_{n}^{(2,k)}\frac{e^{\frac{-2n\mu}{\sigma}}\sigma^{k+1}}{2^{k+1}}\left[\Gamma\left(1+k,\frac{2x}{\sigma}\right)\right]
\right\}
\end{align}
\end{subequations}

The values of coefficients $c_{n+1}$ for $n \geq 1$ are
\begin{equation}
    c_{n+1}=\tilde{f}_{n+1}(0)+\tilde{\tilde{f}}_{n+1}(0)
\end{equation}
where 
\begin{align}\label{massa pos}
&\tilde{f}_{n+1}(0)=\nonumber\\
&=\left\{\frac{1}{(2\sigma)^n}\sum_{k=0}^{n-1}\left[a_{n}^{(1,k)}e^{\frac{2n\mu}{\sigma}}\sigma^{k+1} (-1)^k\left[\left.\Gamma\left(1+k,\frac{-(y-n\mu)}{\sigma}\right)\right|_{0}^{x+n\mu}\right]+b_{n}^{(1,k)}e^{-\frac{2n\mu}{\sigma}}\sigma^{k+1} \left[\left.\Gamma\left(1+k,\frac{(y-n\mu)}{\sigma}\right)\right|_{0}^{x+n\mu}\right]\right]\right.\nonumber\\
&-\frac{e^{\frac{\mu}{\sigma}}}{2^{n+1}\sigma^n}\sum_{k=0}^{n-1}\left[a_{n}^{(1,k)}\frac{e^{\frac{2n\mu}{\sigma}}\sigma^{k+1} (-1)^k}{2^{k+1}}\left[\left.\Gamma\left(1+k,\frac{-2(y-n\mu)}{\sigma}\right)\right|_{0}^{x+n\mu}\right]+b_{n}^{(1,k)}\left.\frac{(y-n\mu)^{k+1}}{k+1}\right|_{0}^{x+n\mu}\right]\nonumber\\
 &+\frac{1}{(2\sigma)^n}\sum_{k=0}^{n-1}\left[a_{n}^{(2,k)}e^{\frac{2n\mu}{\sigma}}\sigma^{k+1} (-1)^k\left[\left.\Gamma\left(1+k,\frac{-(y-n\mu)}{\sigma}\right)\right|_{x+n\mu}^{-\mu}\right]+b_{n}^{(2,k)}e^{-\frac{2n\mu}{\sigma}}\sigma^{k+1} \left[\left.\Gamma\left(1+k,\frac{(y-n\mu)}{\sigma}\right)\right|_{x+n\mu}^{-\mu}\right]\right]\nonumber\\
&-\frac{e^{\frac{\mu}{\sigma}}}{2^{n+1}\sigma^n}\sum_{k=0}^{n-1}\left[a_{n}^{(2,k)}\frac{e^{\frac{2n\mu}{\sigma}}\sigma^{k+1} (-1)^k}{2^{k+1}}\left[\left.\Gamma\left(1+k,\frac{-2(y-n\mu)}{\sigma}\right)\right|_{x+n\mu}^{-\mu}\right]+b_{n}^{(2,k)}\left.\frac{(y-n\mu)^{k+1}}{k+1}\right|_{x+n\mu}^{-\mu}\right]\nonumber\\
&+\frac{e^{-\frac{\mu}{\sigma}}}{2^{n+1}\sigma^n}\sum_{k=0}^{n-1}\left[b_{n}^{(2,k)}\frac{e^{-\frac{2n\mu}{\sigma}}\sigma^{k+1} }{2^{k+1}}\Gamma\left(1+k,-\frac{2(n+1)\mu}{\sigma}\right)\right]+ \left. \left[\left(1-\frac{e^{\frac{\mu}{\sigma}}}{2}\right) c_n\right]\right\} \mathbbm{1}_{[0,-(n+1)\mu)}(x)
 \end{align}
 
and 
\begin{align}\label{massa neg}
 &\tilde{\tilde{f}}_{n+1}(0)=\nonumber\\
 &=\left\{\frac{1}{(2\sigma)^n}\sum_{k=0}^{n-1}\left[a_{n}^{(1,k)}e^{\frac{2n\mu}{\sigma}}\sigma^{k+1} (-1)^k\left[\left.\Gamma\left(1+k,\frac{-(y-n\mu)}{\sigma}\right)\right|_{0}^{-\mu}\right]+b_{n}^{(1,k)}e^{-\frac{2n\mu}{\sigma}}\sigma^{k+1} \left[\left.\Gamma\left(1+k,\frac{(y-n\mu)}{\sigma}\right)\right|_{0}^{-\mu}\right]\right]\right.\nonumber\\
 &-\frac{e^{\frac{\mu}{\sigma}}}{2^{n+1}\sigma^n}\sum_{k=0}^{n-1}\left[a_{n}^{(1,k)}\frac{e^{\frac{2n\mu}{\sigma}}\sigma^{k+1} (-1)^k}{2^{k+1}}\left[\left.\Gamma\left(1+k,\frac{-2(y-n\mu)}{\sigma}\right)\right|_{0}^{-\mu}\right]+b_{n}^{(1,k)}\left.\frac{(y-n\mu)^{k+1}}{k+1}\right|_{0}^{-\mu}\right]\nonumber\\
 &+\frac{e^{-\frac{\mu}{\sigma}}}{2^{n+1}\sigma^n} \sum_{k=0}^{n-1}\left[\left.a_{n}^{(1,k)}\frac{(y-n\mu)^{k+1}}{k+1}\right|_{-\mu}^{x+n\mu}+b_{n}^{(1,k)}\frac{e^{-\frac{2n\mu}{\sigma}}\sigma^{k+1} }{2^{k+1}}\left[\left.\Gamma\left(1+k,\frac{2(y-n\mu)}{\sigma}\right)\right|_{-\mu}^{x+n\mu}\right]\right] \nonumber\\
 &+\frac{e^{-\frac{\mu}{\sigma}}}{2^{n+1}\sigma^n} \sum_{k=0}^{n-1}\left[b_{n}^{(2,k)}\frac{e^{-\frac{2n\mu}{\sigma}}\sigma^{k+1} }{2^{k+1}}\Gamma\left(1+k,\frac{2x}{\sigma}\right)\right]+\left.  \left[\left(1-\frac{e^{\frac{\mu}{\sigma}}}{2}\right)c_n\right]\right\}\mathbbm{1}_{[-(n+1)\mu),\infty)}(x)
 \end{align}

The initial condition for the recursion of the coefficients are 
 \begin{subequations}
\begin{align}
&a_{1}^{(1,0)}=e^{-\frac{\mu+x}{\sigma}}  \\
&b_{1}^{(1,0)}=0\\
&a_{1}^{(2,0)}=0 \label{a_1^(2,0)}\\
&b_{1}^{(2,0)}=e^{\frac{\mu+x}{\sigma}}\label{b_1^(2,0)}\\
&c_1=\frac{e^{-\frac{\mu+x}{\sigma}}}{2}\label{c_1}
\end{align} 
\end{subequations}

\item \underline{If $x+n\mu<0$}, we have $f_n^1(u) \equiv 0$ for each $n$ and the coefficients $a_{n+1}^{(2,j)}$, $b_{n+1}^{(2,j)}$  verify the following recursive relations for $1\leq j\leq n$.
\begin{subequations}\label{coeff a b T2bis}
 \begin{align}
a_{n+1}^{(2,j)}&=0\\
b_{n+1}^{(2,j)}&=e^{\frac{\mu}{\sigma}}\left[\sum_{k=j-1}^{n-1} b_{n}^{(2,k)}\frac{k!}{j!} \left(\frac{\sigma}{2}\right)^{k-j+1}\right]\\
a_{n+1}^{(2,0)}&=0\\
b_{n+1}^{(2,0)}&=e^{\frac{\mu}{\sigma}}\left[\sum_{k=0}^{n-1}
\left[b_{n}^{(2,k)} \left( -\frac{(-n\mu)^{k+1}}{k+1}+k!\left(\frac{\sigma}{2}\right)^{k+1} \right)
\right]\right. \left.-a_{n}^{(2,k)}\frac{e^{\frac{2n\mu}{\sigma}}\sigma^{k+1} (-1)^k}{2^{k+1}}\left[\Gamma\left(1+k,\frac{2n\mu}{\sigma}\right) \right]+(2\sigma)^n f_n^0(0)\right]\\
c_{n+1}&=\frac{1}{(2\sigma)^n}\sum_{k=0}^{n-1}\left[b_{n}^{(2,k)}e^{-\frac{2n\mu}{\sigma}}\sigma^{k+1} \left[\left.\Gamma\left(1+k,\frac{(y-n\mu)}{\sigma}\right)\right|_{0}^{-\mu}\right]\right]-\frac{e^{\frac{\mu}{\sigma}}}{2^{n+1}\sigma^n}\sum_{k=0}^{n-1}\left[b_{n}^{(2,k)}\left.\frac{(y-n\mu)^{k+1}}{k+1}\right|_{0}^{-\mu}\right]\nonumber\\
&+\frac{e^{-\frac{\mu}{\sigma}}}{2^{n+1}\sigma^n}\sum_{k=0}^{n-1}\left[b_{n}^{(2,k)}\frac{e^{-\frac{2n\mu}{\sigma}}\sigma^{k+1} }{2^{k+1}}\Gamma\left(1+k,\frac{-2(n+1)\mu}{\sigma}\right)\right]+  \left[\left(1-\frac{e^{\frac{\mu}{\sigma}}}{2}\right) c_n\right] \label{c_n+1 bis}
\end{align}   
\end{subequations}
with initial conditions given by (\ref{a_1^(2,0)}), (\ref{b_1^(2,0)}) and (\ref{c_1}).

\end{itemize}
\end{thm}

\begin{thm} \label{Theorem Position mu<=-x}
For a Lindley process $\{W_n, n\geq 0\}$ characterized by Laplace increments with location parameter $\mu\leq-x$, the probability density function of the position is given by 
\begin{subequations}\label{fn 3.5}
\begin{align}
f_{n}(u)&=\sum _{i=0}^{1} f_n^i(u), \qquad \qquad u\geq 0, n\geq 1\\
f_{n}(u)&=\delta(u-x)
\end{align} 
\end{subequations}
where 
\begin{subequations}\label{f_n T3}
\begin{align}\label{fn_1 3.5}
&f_n^1(u)=\left\{\sum_{j=0}^{n-1} \frac{1}{(2 \sigma)^n}b_{n}^{(1,j)}(u-(n-1)\mu)^j e^{-\frac{u}{\sigma}}\right\}\mathbbm{1}_{I_n^i}(u) \\
&f_n^0(u)=c_n \delta(u)= c_n\mathbbm{1}_{I_n^0}
\end{align}
\end{subequations}
and $I_n^i$, $i=0,1$, with $\bigcup_{i=0}^{1} I_n^i=[0, +\infty)$, are
 \begin{subequations}
\begin{align}
&I_n^0=\{0\}\\
&I^1_n=(0,+\infty]. 
\end{align}
\end{subequations}

The coefficients  $c_{n+1}$ and $b_{n+1}^{(i,j)}$ verify the following recursive relations for $j=1,\dots, n$
\begin{subequations} \label{abc T3}
\begin{align}
b_{n+1}^{(1,j)}&=e^{\frac{\mu}{\sigma}}\left[\sum_{k=j-1}^{n-1} b_{n}^{((1,k)}\frac{k!}{j!} \left(\frac{\sigma}{2}\right)^{k-j+1}\right]\\
 b_{n+1}^{(i,0)}&=e^{\frac{\mu}{\sigma}}\left[ \sum_{k=0}^{n-1}
\left[b_{n}^{(1,k)} \left(- \frac{[-(n-1)\mu]^{k+1}}{k+1}+k!\left(\frac{\sigma}{2}\right)^{k+1} \right)
\right]+(2\sigma)^{n}f_n^0(0) \right]\\
 c_{n+1}&=-\frac{e^{-\frac{\mu}{\sigma}}}{2^{n+1}\sigma^{n}}\sum_{k=0}^{n-1}b_n^{(1,k)}\frac{e^{-\frac{2(n-1)\mu}{\sigma}}\sigma^{k+1}}{2^{k+1}}\left[\left.\Gamma\left(1+k,\frac{2(y-(n-1)\mu)}{\sigma}\right)\right|^\infty_{-\mu}\right]\nonumber\\
  &-\frac{1}{(2\sigma)^{n}}\sum_{k=0}^{n-1}b_n^{(1,k)}e^{-\frac{(n-1)\mu}{\sigma}}\sigma^{k+1}\left[\left.\Gamma\left(1+k,\frac{(y-(n-1)\mu)}{\sigma}\right)\right|_0^{-\mu}\right]\nonumber\\
  &-\frac{e^{\frac{\mu}{\sigma}}}{2^{n+1}\sigma^{n}}\sum_{k=0}^{n-1}b_n^{(1,k)}\left[\left.\frac{(y-(n-1)\mu)^{k+1}}{k+1}\right|_0^{-\mu}\right]+  \left[\left(1-\frac{e^{\frac{\mu}{\sigma}}}{2}\right) c_n\right].\label{c_n+1 T3}
\end{align}
\end{subequations}

The initial values for the recursion are 
\begin{subequations}
\begin{align}
&c_1=\left(1-\frac{1}{2}e^{\frac{\mu+x}{\sigma}}\right)\\
&b_{1}^{(1,0)}=e^{\frac{\mu+x}{\sigma}}.
\end{align} 
\end{subequations}

\end{thm}
	
\section{First exit time of $W_n$.} \label{Sect. FET}
Let $N_x:=\min_{n>0}\left\lbrace  W_n \geq h|W_0=x \right\rbrace $ be the first exit time (FET) of the Lindley process (\ref{Lindley}) from the domain $[0,h]$ for fixed $h>0$ and let $P(n|x):=\mathbb{P}[N_x=n]$ indicate the probability that the FET is equal to $n$, $n > 0$, given that the process starts in $x\in [0,h)$.

In order to determine the distribution of $N$ computations change according to the sign of $\mu$ and its order with respect to $h$. Theorem \ref{Theorem FPT mu>0 h>mu} gives the distribution for $0<\mu<h$ while Corollary \ref{Theorem FPT mu>0 0<h<mu} considers the case $0<h\leq \mu$. For $\mu<0$, Theorem \ref{Theorem FPT mu<0 } gives the distribution for $-\mu<h$ while Corollary \ref{Theorem FPT mu<0 0<h< -mu} considers the case $0<h<-\mu$ and Theorem \ref{Theorem FPT mu=0} refers to $\mu=0$. 

\begin{thm}\label{Theorem FPT mu>0 h>mu}
	For a Lindley process $\{W_n, n\geq 0\}$ characterized by Laplace increments with location parameter $0<\mu<h$, the probability distribution of the FET through $h$ is given by 
	\begin{equation} \label{Pn}
		P(n|x)=\sum _{i=1}^{\ell_n} P_i(n|x), \qquad \qquad n>0, \quad 0 \leq x < h
	\end{equation}
	where $\ell_n:=\min \{n+1,\min_{r>0}\{r \mu\geq h \}\}$ and 
	\begin{align} 
		&P_i(n|x)=\left\{\sum_{j=0}^{m^i_n-1} \frac{1}{2^n \sigma^{n-1}}\left[\alpha_{n}^{(i,j)}(x+n\mu)^j e^{\frac{x}{\sigma}}+\beta_{n}^{(i,j)}(x+n\mu)^j e^{-\frac{x}{\sigma}}\right]+ \eta_n^{(i,0)}\right\} \mathbbm{1}_{I^i_n}(x)  \qquad 1 \leq i\leq \ell_n \label{Pni}. 
	\end{align}
Here $m_n^i:=\min(n,\ell_n-i+1)$ and we partition the interval $[0,h)$ in
\begin{align}\label{PartitionN}
		&I^i_n=[h-(\ell_n-i+1)\mu,h-(\ell_n-i)\mu) \qquad i=2,\dots , \ell_n\\
		&I^1_n=[0,h-(\ell_n-1)\mu).\nonumber
	\end{align}
	
	The coefficients $\alpha_{n+1}^{(i,j)}$, $\beta_{n+1}^{(i,j)}$ and $\eta_{n+1}^{(i,0)}$ are defined by the following recursive relations for $1 \leq i \leq \ell_{n+1}$ and  $1 \leq j \leq m^i_{n+1}-1$  
\begin{subequations}  
	\begin{align}
&\alpha_{n+1}^{(i,0)}=  e^{\frac{\mu}{\sigma}}\left\{\sum_{j=i^*+1}^{\ell_n}K_n^j\right.+\sum_{k=0}^{m_n^{i^*}-1} 
\left[\alpha_{n}^{(i^*,k)} \left( \frac{(h-(\ell_{n}-i^*)\mu+n\mu)^{k+1}}{k+1} + k!\frac{\sigma^{k+1} (-1)^k}{2^{k+1}}\right)\right.\\ 
&\hspace{1.2cm}-\left.\beta_{n}^{(i^*,k)}\frac{e^{\frac{2n\mu}{\sigma}}\sigma^{k+1} }{2^{k+1}}\left.\Gamma\left(1+k,\frac{2(h-(\ell_{n}-i^*)\mu +n\mu)}{\sigma}\right)\right] -\eta_n^{(i^*,0)} (2 \sigma)^n e^{-\frac{h-(\ell_n-i^*)\mu}{\sigma}}\right\} \nonumber\\
&\alpha_{n+1}^{(i,j)}=-e^{\frac{\mu}{\sigma}}\sum_{k=j-1}^{m_{n}^{i^*}-1}\alpha_{n}^{(i^*,k)}\left(-\frac{\sigma}{2}\right)^{k-j+1}\frac{k!}{j!} \\
&\beta_{n+1}^{(i,0)}=e^{-\frac{\mu}{\sigma}}\left\{\sum_{j=0}^{i^*-1}K_n^j  \right.  +\sum_{k=0}^{m_{n}^{i^*}-1}\left[ -\alpha_{n}^{(i^*,k)}\frac{e^{\frac{-2n\mu}{\sigma}}\sigma^{k+1} (-1)^k}{2^{k+1}}\Gamma\left(1+k,\frac{-2(h-(\ell_{n}-i^*+1)\mu+n\mu)}{\sigma}\right)\right.\\ 
&\hspace{1.2cm}\left.\left.+\beta_{n}^{(i^*,k)}\left(-\frac{(h-(\ell_{n}-i^*+1)\mu+n\mu)^{k+1}}{k+1}+k!\frac{\sigma^{k+1}}{2^{k+1}}\right) \right] -\eta_n^{(i^*,0)} (2 \sigma)^n e^{\frac{h-(\ell_n-i^*+1)\mu}{\sigma}}\right\} \nonumber\\
&\beta_{n+1}^{(i,j)}=e^{-\frac{\mu}{\sigma}}\sum_{k=j-1}^{m_{n}^{i^*}-1}\beta_{n}^{(i^*,k)}\left(\frac{\sigma}{2}\right)^{k-j+1}\frac{k!}{j!}\\
&\eta_{n+1}^{(i,0)}=\eta_{n}^{(i^*,0)} \label{eta},
\end{align}
\end{subequations}
here ${i^*}$ is the index of the interval (\ref{PartitionN}) that contains $x+\mu$, and 
\begin{align}\label{K tot}
 K_n^j=\left\{
 \begin{array}{ll}
       \sum_{k=0}^{m^1_n-1} \left\{ \sigma\left[\alpha_{n}^{(1,k)}(n\mu)^k +\beta_{n}^{(1,k)}(n\mu)^k\right]\right\} + (2\sigma)^n\eta_n^{(1,0)}   \hspace{3cm} j=0\\
       \sum_{k=0}^{m^j_n-1}   \left\{ \alpha_{n}^{(j,k)}\frac{e^{-\frac{2n\mu}{\sigma}}\sigma^{k+1} (-1)^k}{2^{k+1}}\left.\Gamma\left(1+k,\frac{-2(y+n\mu)}{\sigma}\right)\right|_{I_n^j}+\beta_{n}^{(j,k)}\left.\frac{(y+n\mu)^{k+1}}{k+1}\right|_{I_n^{j}}\right\} + \left. \eta_n^{(j,0)}(2\sigma)^n e^{\frac{y}{\sigma}} \right|_{I_n^j} \\
       \hspace{12cm}  j=1,\dots, i^*-1\\
       \sum_{k=0}^{m^j_n-1}   \left\{ \alpha_{n}^{(j,k)}\left.\frac{(y+n\mu)^{k+1}}{k+1}\right|_{I_n^{j}}-\beta_{n}^{(j,k)} \frac{e^{\frac{2n\mu}{\sigma}}\sigma^{k+1} }{2^{k+1}}\left.\Gamma\left(1+k,\frac{2(y+n\mu)}{\sigma}\right)\right|_{I_n^j}\right\} - \left. \eta_n^{(j,0)}(2\sigma)^n e^{-\frac{y}{\sigma}} \right|_{I_n^j}   \\
       \hspace{12cm}  j=i^*+1,\dots, \ell_n       
 \end{array}
  \right.
 \end{align}

The initial values are
\begin{subequations}\label{Ci T FET 1}
\begin{align}
	&\eta_1^{(1,0)}=\beta_{1}^{(1,0)}=\alpha_{1}^{(2,0)}=0 \\
	&\alpha_{1}^{(1,0)}= e^{\frac{\mu-h}{\sigma}}\\
    &\eta_1^{(2,0)}=1\\
	&\beta_{1}^{(2,0)}=-e^{-\frac{\mu-h}{\sigma}}.
\end{align} 
\end{subequations}
\end{thm}

When $\mu\geq h$ the partition reduces to a single interval $[0,h]$ and we get a compact closed form solution.

\begin{corollary}\label{Theorem FPT mu>0 0<h<mu}
    For a Lindley process $\{W_n, n\geq 0\}$ characterized by Laplace increments with location parameter $\mu>0$ and $0<h\leq \mu$, the probability distribution of the FET through $h$ is given by     
    \begin{equation}\label{P mu>0 h<mu}
    P(n|x) = \eta_n+\beta_n e^{-\frac{x}{\sigma}} \qquad \qquad n>0,\quad 0 \leq x < h
    \end{equation}
where
\begin{subequations}
    \begin{align}
    \eta_n &=\delta_{1,n}\\
    \beta_1 &=-\frac{1}{2}e^{\frac{h-\mu}{\sigma}}\\
     \beta_{n}&=\frac{e^{\frac{h-(n-1)\mu}{\sigma}}}{2}\left(\frac{1}{2}+\frac{h}{2\sigma}\right)^{n-2}-\frac{e^{\frac{h-n\mu}{\sigma}}}{2}\left(\frac{1}{2}+\frac{h}{2\sigma}\right)^{n-1}, \qquad n \geq 2.\label{beta2 Teor4.2}
    \end{align}
    \end{subequations}
\end{corollary}

\begin{thm}  \label{Theorem FPT mu<0 }
    For a Lindley process $\{W_n, n\geq 0\}$ characterized by Laplace increments with location parameter $\mu<0$, the probability distribution of the FET through $h$ is given by  
  \begin{equation}\label{P_n mu<0}
    P(n|x)=\sum _{i=1}^{h_n} P_i(n|x), \qquad \qquad n>0, \quad 0 \leq x < h
  \end{equation}
  where $h_n:=\min\{n,\min_{r>0}\{r(-\mu)>h\}\}$ and, for  $1 \leq i\leq h_n  $, 
  \begin{align} \label{P_ni mu<0}
    &P_i(n|x)=\left\{\sum_{j=0}^{n-1} \frac{1}{2^n \sigma^{n-1}}\left[\alpha_{n}^{(i,j)}(x+(n-1)\mu)^j e^{\frac{x}{\sigma}}+\beta_{n}^{(i,j)}(x+(n-1)\mu)^j e^{-\frac{x}{\sigma}}\right]+ \eta_n^{i} \right\}\mathbbm{1}_{I^i_n}(x).
  \end{align}
Here we partition the interval $[0,h)$ as
\begin{align}\label{PartitionN1}
		&I^i_n=[-(i-1)\mu, -i\mu) \qquad i=1,...,h_n-1\\
		&I^{h_n}_n=[-(h_n-1)\mu,h).\nonumber
	\end{align}
The coefficients $\alpha_{n+1}^{(i,j)}$, $\beta_{n+1}^{(i,j)}$ and $\eta_{n+1}^{i}$ are defined by the recursive relations, for $1<i\leq h_{n+1}$ and $0<j\leq n$
\begin{subequations}\label{coeff N4.4}
 \begin{align}
&\alpha_{n+1}^{(i,0)}=  e^{\frac{\mu}{\sigma}}\left\{  \sum_{j=i^*+1}^{h_n}K_n^j\right. +\sum_{k=0}^{n-1} 
\left[\alpha_{n}^{(i^*,k)} \left( \frac{((-i^*+n-1)\mu)^{k+1}}{k+1} + k!\frac{\sigma^{k+1} (-1)^k}{2^{k+1}}\right)\right.\\ 
&\hspace{3.2cm}-\left.\beta_{n}^{(i^*,k)}\frac{e^{\frac{2(n-1)\mu}{\sigma}}\sigma^{k+1} }{2^{k+1}}\left.\Gamma\left(1+k,\frac{2(-i^*+n-1)\mu)}{\sigma}\right)\right] -\eta^{i^*} (2 \sigma)^n e^{\frac{i^*\mu}{\sigma}}\right\}\nonumber\\
&\alpha_{n+1}^{(i,j)}=-e^{\frac{\mu}{\sigma}}\sum_{k=j-1}^{n-1}\alpha_{n}^{(i^*,k)}\left(-\frac{\sigma}{2}\right)^{k-j+1}\frac{k!}{j!} \\
&\beta_{n+1}^{(i,0)}=e^{-\frac{\mu}{\sigma}}\left\{\sum_{j=0}^{i^*-1}K_n^j  \right. +\sum_{k=0}^{n-1}\left[ -\alpha_{n}^{(i^*,k)}\frac{e^{\frac{-2(n-1)\mu}{\sigma}}\sigma^{k+1} (-1)^k}{2^{k+1}}\Gamma\left(1+k,\frac{-2(-i^*+1+n-1)\mu)}{\sigma}\right)\right.\\ 
&\left.
\hspace{3.4cm}\left.+\beta_{n}^{(i^*,k)}\left(-\frac{((-i^*+1+n-1)\mu)^{k+1}}{k+1}+k!\frac{\sigma^{k+1}}{2^{k+1}}\right) \right] -\eta^{i^*} (2 \sigma)^n e^{-\frac{(i^*-1)\mu}{\sigma}} \right\}\nonumber\\
&\beta_{n+1}^{(i,j)}=e^{-\frac{\mu}{\sigma}}\sum_{k=j-1}^{n-1}\beta_{n}^{(i^*,k)}\left(\frac{\sigma}{2}\right)^{k-j+1}\frac{k!}{j!} \\
&\eta_{n+1}^{i}=\eta_{n}^{i^*}
\end{align}
\end{subequations}
where $I_n^{i^*}$ is the interval that contains $\min(x+\mu,h)$.\\
For $i=1$ the coefficients are defined by the recursive relations
\begin{subequations}\label{coeff N4.4 i=1} 
\begin{align}
&\alpha_{n+1}^{(1,0)}=e^{\frac{\mu}{\sigma}}\left\{ \sum_{j=1}^{h_n}K_n^j -K_n^0 \right\}
\end{align}
\begin{align}
&\alpha_{n+1}^{(1,j)}=0 \qquad j=1\dots n\\
&\beta_{n+1}^{(1,j)}=0 \qquad j=0\dots n\\
&\eta_{n+1}^{1}=\frac{K_n^0}{2^n\sigma^n}
\end{align}
\end{subequations}
where
\begin{align}\label{K tot mu<0}
 K_n^j=\left\{
 \begin{array}{ll} 
\sum_{k=0}^{n-1} \left\{ \sigma\left[\alpha_{n}^{(1,k)}((n-1)\mu)^k +\beta_{n}^{(1,k)}((n-1)\mu)^k\right]\right\} + (2\sigma)^n\eta_n^{1}
 \hspace{2.9cm} j=0\\
        \sum_{k=0}^{n-1}   \left\{ \alpha_{n}^{(j,k)}\frac{e^{-\frac{2n\mu}{\sigma}}\sigma^{k+1} (-1)^k}{2^{k+1}}\left.\Gamma\left(1+k,\frac{-2(y+(n-1)\mu)}{\sigma}\right)\right|_{I_n^j}+\beta_{n}^{(j,k)}\left.\frac{(y+(n-1)\mu)^{k+1}}{k+1}\right|_{I_n^{j}}\right\} + \left. \eta_n^{j}(2\sigma)^n e^{\frac{y}{\sigma}} \right|_{I_n^j} \\
       \hspace{12.5cm}  j=1,\dots, i^*-1\\
       \sum_{k=0}^{n-1}   \left\{ \alpha_{n}^{(j,k)}\left.\frac{(y+(n-1)\mu)^{k+1}}{k+1}\right|_{I_n^{j}}-\beta_{n}^{(j,k)} \frac{e^{\frac{2n\mu}{\sigma}}\sigma^{k+1} }{2^{k+1}}\left.\Gamma\left(1+k,\frac{2(y+(n-1)\mu)}{\sigma}\right)\right|_{I_n^j}\right\} - \left. \eta_n^{j}(2\sigma)^n e^{-\frac{y}{\sigma}} \right|_{I_n^j}   \\
       \hspace{12.5cm}  j=i^*+1,\dots, h_n       
 \end{array}
  \right.
 \end{align}

The coefficients for the base case of the recursion are given by
\begin{subequations}  
\begin{align}
&\alpha_{1}^{(1,0)}=e^{\frac{\mu-h}{\sigma}}\label{alpha_1 T3}\\
&\beta_{1}^{(1,0)}=\eta_{1}^{1}=0
\end{align}
\end{subequations}

\end{thm}

When $-\mu\geq h$ in the previous theorem we get that the partition is formed of a single interval $[0,h]$. This simplifies a lot the computations and we get the following:
 
\begin{corollary} \label{Theorem FPT mu<0 0<h< -mu}
    For a Lindley process $\{W_n, n\geq 0\}$ characterized by Laplace increments with location parameter $\mu<0$, and $0<h\leq -\mu$, the probability distribution of the FET through $h$ is given by     
    \begin{equation}\label{P_n mu<0 h<-mu}
    P(n|x) = \eta_n+\alpha_n e^{\frac{x}{\sigma}}   \qquad \qquad n>0, \quad 0 \leq x < h
    \end{equation}
where, for $n>1$
\begin{align}
\alpha_{n+1}&= \frac{e^{\frac{\mu}{\sigma}}}{2}\left(\frac{h}{\sigma}-1\right)\alpha_n-\frac{e^{\frac{\mu-h}{\sigma}}}{2}\eta_n \label{alpha T4}\\ 
 \eta_{n+1}&=\alpha_n+\eta_n\label{eta T4}
\end{align}
and the initial values are
\begin{align}
\alpha_{1}&= \frac{1}{2}e^{\frac{\mu-h}{\sigma}} \\ 
\eta_{1}&=0
\end{align}
\end{corollary}

\begin{thm}\label{Theorem FPT mu=0} 
For a Lindley process $\{W_n, n\geq 0\}$ characterized by Laplace increments with location parameter $\mu=0$, the probability distribution of the FET through $h$ is given by 
\begin{align} 
	&P(n|x)=\frac{1}{2^n \sigma^{n-1}}\left[\sum_{j=0}^{n-1} \alpha_{n}^{(1,j)}x^j e^{\frac{x}{\sigma}}+\sum_{j=0}^{n-2} \beta_{n}^{(1,j)}x^j e^{-\frac{x}{\sigma}}\right] \qquad \qquad n>0, \quad 0 \leq x < h
\end{align}
The coefficients satisfy the following recursive relations  for $n>1$
\begin{subequations} \label{alfa beta mu=0}
\begin{align}
    \alpha_{n+1}^{(1,j)}&=-\sum_{k=j-1}^{n-1}\alpha_n^{(1,k)}\left(-\frac{\sigma}{2}\right)^{k-j+1}\frac{k!}{j!} \qquad\qquad j=1, \dots, n\\
    \beta_{n+1}^{(1,j)}&=\sum_{k=j-1}^{n-2}\beta_n^{(1,k)}\left(\frac{\sigma}{2}\right)^{k-j+1}\frac{k!}{j!} \qquad\qquad j=1, \dots, n-1
    \end{align}
    \begin{align}
    \alpha_{n+1}^{(1,0)}&=\sum_{k=0}^{n-1}\alpha_n^{(1,k)}\left(\frac{h^{k+1}}{k+1}+\frac{\sigma^{k+1}}{2^{k+1}}(-1)^kk!\right)+\sum_{k=0}^{n-2}\beta_{n}^{(1,k)}\frac{\sigma^{k+1}}{2^{k+1}}\Gamma\left(1+k,\frac{2h}{\sigma}\right)\\
    \beta_{n+1}^{(1,0)}&=-\sum_{k=0}^{n-1}\alpha_n^{(1,k)}\frac{\sigma^{k+1}}{2^{k+1}}(-1)^kk!+\sum_{k=0}^{n-2}\beta_n^{(1,k)}\frac{\sigma^{k+1}}{2^{k+1}}k!+\sigma(\alpha_n^{(1,0)}+\beta_n^{(1,0)})
\end{align}
\end{subequations}

with initial conditions
\begin{subequations}
    \begin{align}
\alpha_1^{(1,0)}&=e^{-\frac{h}{\sigma}}\\
\beta_1^{(1,0)}&=0
\end{align}
\end{subequations}

\end{thm}

\begin{rem}
In Theorem \ref{Theorem FPT mu>0 h>mu} and Theorem \ref{Theorem FPT mu<0 }, when $\mu$ is not small, the value of $i$ such that $i|\mu|>h$ is small.  Hence the sums include very few terms since $l_n$ and  $k_n$ coincide with $i$ soon.
\end{rem}
\begin{rem}
  Note that from Corollary \ref{Theorem FPT mu<0 0<h< -mu} we easily observe the exponential decay of the tails of the FET probability function.  Since larger values of $\mu$ facilitate the crossing of the boundary, this tail result holds for any choice of $\mu$.
\end{rem}

\section{Proofs of Theorems on the position}

\subsection{Proof of Lemma \ref{Lemma1}}
\begin{proof}
   The thesis follows immediately determining the distribution of $x+Z_1$ and applying  definition (\ref{Lindley}).
\end{proof}

\subsection{Proof of Theorem \ref{Theorem Position mu>=0}}

\begin{proof}
We proceed by induction. 
 
 {\bf Case n=1}
 Since $\mu>0$, we have $x>-\mu$. Hence, this case follows from Lemma \ref{Lemma1} recognising the partition of $[0,+\infty)$ as $\bigcup_{i=0}^{2} I_n^i$ where 
\begin{align}
&I_1^0=\{0\}\\
&I_1^1=(0,\mu+x]\\
&I_1^2=(\mu+x,+\infty)
\end{align} 
and rewriting the corresponding density function (\ref{f_1 generale}) as
\begin{align}\label{f_1}
f_1(u)&=f_{W_1}(u)\\
&=\left\lbrace 
\begin{array}{ll}
\frac{e^{-\frac{\mu+x}{\sigma}}}{2}  &u=0 \nonumber\\
\frac{1}{2\sigma}e^{-\frac{\mu-u+x}{\sigma}}  &0<u \leq x+\mu \nonumber\\
\frac{1}{2\sigma}e^{\frac{\mu-u+x}{\sigma}} & u>x+\mu
\end{array}
\right. \nonumber\\
&=\frac{e^{-\frac{\mu+x}{\sigma}}}{2}  \mathbbm{1}_{I_1^0}(u) +\frac{1}{2\sigma}e^{-\frac{\mu-u+x}{\sigma}} \mathbbm{1}_{I_1^1}(u)+
\frac{1}{2\sigma}e^{\frac{\mu-u+x}{\sigma}} \mathbbm{1}_{I_1^2}(u) \nonumber\\
&=\sum _{i=0}^{2} f_1^i(u), \qquad \qquad u\geq 0 \nonumber
\end{align} 
where $\mathbbm{1}_{I_1^0}(u)=\delta(u)$. 
The proof is completed recognizing the initial conditions (\ref{CI mu>0})

 {\bf Case n}
 
 We assume that (\ref{fn}) holds for $n$, and we show that it holds for $n+1$, for $n\geq 1$. 
 
Conditioning on the position reached at time $n$ we can write  the probability distribution function of the position at time $n+1$ as
 \begin{equation}
 F_{n+1}(u)=\mathbb{P}(W_{n+1} \leq u)=\int_0^{+\infty} \mathbb{P}(W_{n+1} \leq u|W_{n}=y) \mathbb{P}(W_n\in dy)  \qquad n=1,2, \dots
 \end{equation}
Differentiating with respect to $u$ we get the equivalent relation for the probability density function of the position at time $n+1$
 \begin{equation}\label{CK1}
 f_{n+1}(u)=\frac{\partial}{\partial u}\mathbb{P}(W_{n+1} \leq u)=\int_0^{+\infty} f(u,n+1|y,n) f_n(y) dy, 	\qquad n=1,2,\dots
 \end{equation}
 Using (\ref{denstransiz}) and (\ref{fn}), equation (\ref{CK1}) becomes
 \begin{align}
 f_{n+1}(u)&=\sum _{i=0}^{n+1} \left[\int_0^{+\infty} \frac{e^{-\frac{|u-y-\mu|}{\sigma}}}{2\sigma} f_n^i(y) dy\right]\mathbbm{1}_{(0,+\infty)}(u)+ \sum _{i=0}^{n+1} \left[\int_0^{+\infty} \frac{e^{-\frac{y+\mu}{\sigma}}}{2} f_n^i(y) dy\right]\delta(u) \label{Continua}
 \end{align}

Note that $f_{n+1}(u)$ is decomposed in a continuous part and in a discrete one.

  Let us observe that when we open the modulus in (\ref{Continua}) we add a further interval in the partition $\{I_{n}^j$, $j=0, \dots , n+1\}$. Indeed, since $(i-2)\mu<u-\mu<(i-1) \mu$ implies   $(i-1)\mu<u< i\mu$, it induces a shift in the partition intervals (\ref{I_n}) and we get  $\{I_{n+1}^j$, $j=0, \dots , n+2\}$. Hence $f_{n+1}(u)$ can be expressed in terms of $ f_{n+1}^i(u)$, according to the new partition.
\begin{equation}
f_{n+1}(u)=\sum _{i=0}^{n+2} f_{n+1}^i(u), \qquad \qquad u\geq 0
\end{equation}
where 
\begin{align}
&f_{n+1}^i(u)=\left\{\sum _{j=0}^{n+1} \left[\int_0^{+\infty} \frac{e^{-\frac{|u-y-\mu|}{\sigma}}}{2\sigma} f_n^j(y) dy\right]\right\} \mathbbm{1}_{I_{n+1}^i}(u), \qquad i=1,\dots n+2 \label{cont}\\
&f_{n+1}^0(u)=\left\{\sum _{j=0}^{n+1} \left[ \int_0^{+\infty} \frac{e^{-\frac{y+\mu}{\sigma}}}{2} f_n^j(y) dy\right]\right\}\delta(u)=c_{n+1} \delta(u)\label{discrete}
\end{align}

Let us consider (\ref{cont}). It holds
\begin{align}
f_{n+1}^1(u)&= \left\{\frac{e^{\frac{(u-\mu)}{\sigma}}}{2\sigma} c_n+
\sum _{j=1}^{n+1} \int_{I_{n}^j} \frac{e^{\frac{(u-y-\mu)}{\sigma}}}{2\sigma} f_n^j(y) dy\right \} \mathbbm{1}_{I_{n+1}^1}(u) \label{f^1_n+1}\\
f_{n+1}^i(u)&= \left\{\frac{e^{-\frac{(u-\mu)}{\sigma}}}{2\sigma} c_n+\sum _{j=1}^{i-2} \int_{I_{n}^j} \frac{e^{-\frac{(u-y-\mu)}{\sigma}}}{2\sigma} f_n^j(y) dy+
\int_{(i-2)\mu}^{u-\mu} \frac{e^{-\frac{(u-y-\mu)}{\sigma}}}{2\sigma} f_n^{i-1}(y) dy\right.  \label{f^i_n+1}\\
 &\qquad\left.+\int_{u-\mu}^{(i-1)\mu}  \frac{e^{\frac{(u-y-\mu)}{\sigma}}}{2\sigma} f_n^{i-1}(y) dy+
\sum _{j=i}^{n+1} \int_{I_{n}^j} \frac{e^{\frac{(u-y-\mu)}{\sigma}}}{2\sigma} f_n^j(y) dy \right \} \mathbbm{1}_{I_{n+1}^i}(u) \qquad i=2,\dots,n \nonumber \\
f_{n+1}^{n+1}(u)&= \left\{\frac{e^{-\frac{(u-\mu)}{\sigma}}}{2\sigma} c_n +\sum _{j=1}^{n-1} \int_{I_{n}^j} \frac{e^{-\frac{(u-y-\mu)}{\sigma}}}{2\sigma} f_n^j(y) dy+
\int_{(n-1)\mu}^{u-\mu} \frac{e^{-\frac{(u-y-\mu)}{\sigma}}}{2\sigma} f_n^{n}(y) dy\right.\label{f^n+1_n+1}\\
& \qquad\left.+\int_{u-\mu}^{x+n\mu} \frac{e^{\frac{(u-y-\mu)}{\sigma}}}{2\sigma} f_n^{n}(y) dy+
\int_{I_{n}^{n+1}} \frac{e^{\frac{(u-y-\mu)}{\sigma}}}{2\sigma} f_n^{n+1}(y) dy \right \} \mathbbm{1}_{I_{n+1}^{n+1}}(u)\nonumber
\end{align}
\begin{align}
f_{n+1}^{n+2}(u)&= \left\{\frac{e^{-\frac{(u-\mu)}{\sigma}}}{2\sigma} c_n +\sum _{j=1}^{n} \int_{I_{n}^j} \frac{e^{-\frac{(u-y-\mu)}{\sigma}}}{2\sigma} f_n^j(y) dy+
\int_{x+n\mu}^{u-\mu} \frac{e^{-\frac{(u-y-\mu)}{\sigma}}}{2\sigma} f_n^{n}(y) dy\right.\label{f^n+2_n+1}\\
& \qquad \left. +\int_{u-\mu}^{+\infty} \frac{e^{\frac{(u-y-\mu)}{\sigma}}}{2\sigma} f_n^{n}(y) dy \right\} \mathbbm{1}_{I_{n+1}^{n+1}}(u)\nonumber
\end{align}

Let us now focus on (\ref{f^i_n+1}) that holds for $2 \leq i \leq n$. Using the inductive property on $f_n^j(y)$, i.e. substituting (\ref{fniPos}) in (\ref{f^i_n+1}), we get
\begin{subequations}
 \begin{align}
f_{n+1}^i(u)&= \left\{\frac{e^{\frac{-(u-\mu)}{\sigma}}}{2\sigma} c_n\right.\\
&+\sum _{j=1}^{i-2} \frac{1}{(2 \sigma)^n}\sum_{k=0}^{m^j_n-1}  \int_{I_{n}^j} \frac{e^{-\frac{(u-y-\mu)}{\sigma}}}{2\sigma} \left[a_{n}^{(j,k)}(y-n\mu)^k e^{\frac{y}{\sigma}}+b_{n}^{(j,k)}(y-n\mu)^k e^{-\frac{y}{\sigma}}\right]   dy\\
&+\frac{1}{(2 \sigma)^n}\sum_{k=0}^{m^{i-1}_n-1} \int_{(i-2)\mu}^{u-\mu} \frac{e^{-\frac{(u-y-\mu)}{\sigma}}}{2\sigma} \left[a_{n}^{((i-1),k)}(y-n\mu)^k e^{\frac{y}{\sigma}}+b_{n}^{((i-1),k)}(y-n\mu)^k e^{-\frac{y}{\sigma}}\right]dy\\
&+\frac{1}{(2 \sigma)^n}\sum_{k=0}^{m^{i-1}_n-1} \int_{u-\mu}^{(i-1)\mu} \frac{e^{\frac{(u-y-\mu)}{\sigma}}}{2\sigma} \left[a_{n}^{((i-1),k)}(y-n\mu)^k e^{\frac{y}{\sigma}}+b_{n}^{((i-1),k)}(y-n\mu)^k e^{-\frac{y}{\sigma}}\right] dy\\
&\left.+\sum _{j=i}^{n+1} \frac{1}{(2 \sigma)^n}\sum_{k=0}^{m^{j}_n-1} \int_{I_{n}^j} \frac{e^{\frac{(u-y-\mu)}{\sigma}}}{2\sigma} \left[a_{n}^{(j,k)}(y-n\mu)^k e^{\frac{y}{\sigma}}+b_{n}^{(j,k)}(y-n\mu)^k e^{-\frac{y}{\sigma}}\right] dy \right \} \mathbbm{1}_{I_{n+1}^i}(u)
\end{align}
\end{subequations}
which can be written as
\begin{subequations}
 \begin{align}
f_{n+1}^i(u)&=\left\{ \frac{e^{-\frac{(u-\mu)}{\sigma}}}{2\sigma} c_n \right.\\
&+\sum _{j=1}^{i-2} \frac{e^{-\frac{(u-\mu)}{\sigma}}}{(2 \sigma)^{n+1}}\sum_{k=0}^{m^j_n-1}  \int_{I_{n}^j} \left[a_{n}^{(j,k)}(y-n\mu)^k e^{\frac{2y}{\sigma}}+b_{n}^{(j,k)}(y-n\mu)^k \right] dy\\
&+\frac{e^{-\frac{(u-\mu)}{\sigma}}}{(2 \sigma)^{n+1}}\sum_{k=0}^{m^{i-1}_n-1} \int_{(i-2)\mu}^{u-\mu} \left[a_{n}^{((i-1),k)}(y-n\mu)^k e^{\frac{2y}{\sigma}}+b_{n}^{((i-1),k)}(y-n\mu)^k \right]dy\\
&+\frac{e^{\frac{(u-\mu)}{\sigma}}}{(2 \sigma)^{n+1}}\sum_{k=0}^{m^{i-1}_n-1} \int_{u-\mu}^{(i-1)\mu}  \left[a_{n}^{((i-1),k)}(y-n\mu)^k +b_{n}^{((i-1),k)}(y-n\mu)^k e^{-\frac{2y}{\sigma}}\right] dy\\
&\left.+\sum _{j=i}^{n+1} \frac{e^{\frac{(u-\mu)}{\sigma}}}{(2 \sigma)^{n+1}}\sum_{k=0}^{m^j_n-1} \int_{I_{n}^j} \left[a_{n}^{(j,k)}(y-n\mu)^k +b_{n}^{(j,k)}(y-n\mu)^k e^{-\frac{2y}{\sigma}}\right] dy \right \} \mathbbm{1}_{I_{n+1}^i}(u)
\end{align}
\end{subequations}

Computing the integrals we get
\begin{subequations}
 \begin{align}
f_{n+1}^i(u)&= \left\{\frac{e^{-\frac{(u-\mu)}{\sigma}}}{2\sigma} c_n \right.\\
&+\sum _{j=1}^{i-2} \frac{e^{-\frac{(u-\mu)}{\sigma}}}{(2 \sigma)^{n+1}}\sum_{k=0}^{m^j_n-1}   \left\{a_{n}^{(j,k)}\frac{e^{\frac{2n\mu}{\sigma}}\sigma^{k+1} (-1)^k}{2^{k+1}}\left[\left.\Gamma\left(1+k,\frac{-2(y-n\mu)}{\sigma}\right)\right|_{I_{n}^j} \right]+b_{n}^{(j,k)}\left.\frac{(y-n\mu)^{k+1}}{k+1}\right|_{I_{n}^j} \right\}
 \end{align}
    \begin{align}
&+\frac{e^{-\frac{(u-\mu)}{\sigma}}}{(2 \sigma)^{n+1}}\sum_{k=0}^{m^{i-1}_n-1} \left\{a_{n}^{((i-1),k)}
\frac{e^{\frac{2n\mu}{\sigma}}\sigma^{k+1} (-1)^k}{2^{k+1}}\left[\left.\Gamma\left(1+k,\frac{-2(y-n\mu)}{\sigma}\right)\right|_{(i-2) \mu}^{u-\mu}\right]
+b_{n}^{((i-1),k)}\left.\frac{(y-n\mu)^{k+1}}{k+1}\right|_{(i-2)\mu}^{u-\mu} \right\}\\
&+\frac{e^{\frac{(u-\mu)}{\sigma}}}{(2 \sigma)^{n+1}}\sum_{k=0}^{m^{i-1}_n-1}  \left\{a_{n}^{((i-1),k)}\left.\frac{(y-n\mu)^{k+1}}{k+1}\right|_{u-\mu}^{(i-1)\mu} -b_{n}^{((i-1),k)}\frac{e^{\frac{-2n\mu}{\sigma}}\sigma^{k+1}}{2^{k+1}}\left[\left.\Gamma\left(1+k,\frac{2(y-n\mu)}{\sigma}\right)\right|_{u-\mu}^{(i-1)\mu}\right]
\right\}\\
&\left.+\sum _{j=i}^{n+1} \frac{e^{\frac{(u-\mu)}{\sigma}}}{(2 \sigma)^{n+1}}\sum_{k=0}^{m^j_n-1} \left\{a_{n}^{(j,k)}\left.\frac{(y-n\mu)^{k+1}}{k+1}\right|_{I_{n}^j} -b_{n}^{(j,k)} \frac{e^{\frac{-2n\mu}{\sigma}}\sigma^{k+1}}{2^{k+1}}\left[\left.\Gamma\left(1+k,\frac{2(y-n\mu)}{\sigma}\right)\right|_{I_{n}^j}\right]
\right\} \right \} \mathbbm{1}_{I_{n+1}^i}(u)\label{converg}
\end{align}
\end{subequations}
where $\Gamma(n,y)$ denotes the incomplete Gamma function  
\begin{equation}
    \Gamma(n,y)=\int_y^\infty t^{n-1} e^{-t} dt.
\end{equation}
and $\left.\Gamma(n,\phi(y))\right|_{[a,b]}=\Gamma(n,\phi(b))-\Gamma(n,\phi(a))$.

The divergence of $(y-n\mu)^{k+1}$ on the last interval $I_n^{n+1}=(n\mu+x, +\infty)$ in (\ref{converg}) will not be a problem because from the recurrent expression of the coefficients we will get $a_n^{(n+1,k)}\equiv 0$ for each admissible $k$. 

Using (\ref{A_j}) and the expansion \cite{DLMF}
\begin{equation} \label{gammaExpansion}
    \Gamma(n+1,x)=n! e^{-x} \sum_{r=0}^n \frac{x^r}{r!}
\end{equation}
we get 

\begin{subequations}
 \begin{align}
f_{n+1}^i(u)&= \left\{\frac{e^{-\frac{(u-\mu)}{\sigma}}}{2\sigma} c_n \right.\\
&+\frac{e^{-\frac{(u-\mu)}{\sigma}}}{(2 \sigma)^{n+1}}\sum _{j=1}^{i-2}  A^{(1)}_j\\
&+\frac{e^{-\frac{(u-\mu)}{\sigma}}}{(2 \sigma)^{n+1}}\sum_{k=0}^{m^{i-1}_n-1} \left\{a_{n}^{((i-1),k)}
e^{\frac{2y}{\sigma}}k! \left.
 \sum_{r=0}^k \frac{(-1)^{k+r}(2)^{r-k-1}}{r! \sigma^{r-k-1}}(y-n\mu)^r
\right|_{(i-2) \mu}^{u-\mu}
+b_{n}^{((i-1),k)}\left.\frac{(y-n\mu)^{k+1}}{k+1}\right|_{(i-2)\mu}^{u-\mu} \right\}\\
&+\frac{e^{\frac{(u-\mu)}{\sigma}}}{(2 \sigma)^{n+1}}\sum_{k=0}^{m^{i-1}_n-1}  \left\{a_{n}^{((i-1),k)}\left.\frac{(y-n\mu)^{k+1}}{k+1}\right|_{u-\mu}^{(i-1)\mu} +b_{n}^{((i-1),k)} e^{\frac{-2y}{\sigma}}k!
\left.
\sum_{r=0}^k \frac{-(2)^{r-k-1}}{r! \sigma^{r-k-1}}(y-n\mu)^r
\right|_{u-\mu}^{(i-1)\mu}
\right\}\\
&\left.+\frac{e^{\frac{(u-\mu)}{\sigma}}}{(2 \sigma)^{n+1}}\sum _{j=i}^{n+1}  A^{(2)}_j  \right \} \mathbbm{1}_{I_{n+1}^i}(u) \hspace{6cm}
\end{align}
\end{subequations}

Collecting powers $\left(u-(n+1)\mu\right)^j e^{\frac{u}{\sigma}}$ and $(u-(n+1)\mu)^j e^{-\frac{u}{\sigma}}$, we get for $i>1$ and $j=1,\dots, m_{n+1}^i$

\begin{subequations}
 \begin{align}
a_{n+1}^{(i,j)}&=e^{-\frac{\mu}{\sigma}}\left[-\frac{a_{n}^{((i-1),(j-1))}}{j}+\sum_{k=j}^{m^{i-1}_n-1}(-1)^{k+j} a_{n}^{((i-1),k)}
\frac{k!}{j!}   \left(\frac{2}{\sigma}\right)^{j-k-1}
\right]\\
&=e^{-\frac{\mu}{\sigma}}\left[\sum_{k=j-1}^{m^{i-1}_n-1} (-1)^{k+j} a_{n}^{((i-1),k)}
\frac{k!}{j!}   \left(\frac{2}{\sigma}\right)^{j-k-1}
\right]\\
b_{n+1}^{(i,j)}&=e^{\frac{\mu}{\sigma}}\left[\frac{b_{n}^{((i-1),(j-1))}}{j}+\sum_{k=j}^{m^{i-1}_n-1} b_{n}^{((i-1),k)}\frac{k!}{j!} \left(\frac{2}{ \sigma}\right)^{j-k-1}\right]\\
&=e^{\frac{\mu}{\sigma}}\left[\sum_{k=j-1}^{m^{i-1}_n-1} b_{n}^{((i-1),k)}\frac{k!}{j!} \left(\frac{2}{ \sigma}\right)^{j-k-1}\right]
\end{align}
\end{subequations}
while the coefficients of $e^{\frac{u}{\sigma}}$ and $e^{-\frac{u}{\sigma}}$ are respectively

\begin{subequations}
 \begin{align}
a_{n+1}^{(i,0)}&=e^{-\frac{\mu}{\sigma}}\left[\sum _{j=i}^{n+1}  A^{(2)}_j +\sum_{k=0}^{m^{i-1}_n-1}\left[a_{n}^{((i-1),k)}\left(  \frac{((i-1-n)\mu)^{k+1}}{k+1}+(-1)^{k}k!\left(\frac{\sigma}{2}\right)^{k+1} \right)\right.\right.\\
&\left.-b_{n}^{((i-1),k)}\frac{e^{\frac{-2n\mu}{\sigma}}\sigma^{k+1}}{2^{k+1}}\left[\Gamma\left(1+k,\frac{2((i-1-n)\mu)}{\sigma}\right)\right]\right]\\
b_{n+1}^{(i,0)}&=e^{\frac{\mu}{\sigma}}\left[\sum _{j=1}^{i-2}  A^{(1)}_j- \sum_{k=0}^{m^{i-1}_n-1}
\left[b_{n}^{((i-1),k)} \left( \frac{((i-2-n)\mu)^{k+1}}{k+1}-k!\left(\frac{\sigma}{2}\right)^{k+1} \right)
\right]\right.\\
& \left.-a_{n}^{((i-1),k)}\frac{e^{\frac{2n\mu}{\sigma}}\sigma^{k+1} (-1)^k}{2^{k+1}}\left[\Gamma\left(1+k,\frac{-2((i-1-n)\mu)}{\sigma}\right) \right]+c_n\right]
\end{align}
\end{subequations}
In a similar way, expanding (\ref{f^1_n+1}), (\ref{f^n+1_n+1}) and (\ref{f^n+2_n+1}), we get the thesis. 

In order to find the recursive expression for $c_n$,  let us consider (\ref{discrete})
\begin{subequations}
 \begin{align}
    c_{n+1}&=\int_0^{+\infty} \frac{e^{-\frac{y+\mu}{\sigma}}}{2} f_n(y) dy\\
        &=\sum _{j=0}^{n+1} \left[ \int_0^{+\infty} \frac{e^{-\frac{y+\mu}{\sigma}}}{2} f_n^j(y) dy\right]\\
    &=\frac{e^{-\frac{\mu}{\sigma}}}{2}\left\{\sum _{j=1}^{n+1} \sum_{k=0}^{m^i_n-1} \left[ \int_{I_n^j} \frac{e^{-\frac{y}{\sigma}}}{(2 \sigma)^n} \left[a_{n}^{(j,k)}(y-n\mu)^k e^{\frac{y}{\sigma}}+b_{n}^{(j,k)}(y-n\mu)^k e^{-\frac{y}{\sigma}}\right] dy\right]+c_n\right\}\\
    &=\frac{e^{-\frac{\mu}{\sigma}}}{2}\left\{\sum _{j=1}^{n+1} \sum_{k=0}^{m^i_n-1} \left[ \int_{I_n^j} \frac{1}{(2 \sigma)^n} \left[a_{n}^{(j,k)}(y-n\mu)^k +b_{n}^{(j,k)}(y-n\mu)^k e^{-\frac{2y}{\sigma}}\right] dy\right]+c_n\right\}\\
    &=\frac{e^{-\frac{\mu}{\sigma}}}{2(2 \sigma)^n}\left\{\sum _{j=1}^{n+1} \sum_{k=0}^{m^i_n-1}  \left[\left.a_{n}^{(j,k)}\frac{(y-n\mu)^{k+1}}{k+1}\right|_{I_n^j} -b_{n}^{(j,k)}\frac{e^{-\frac{2n\mu}{\sigma}}\sigma^{k+1}}{2^{k+1}}\left.\Gamma\left(1+k,\frac{2(y-n\mu)}{\sigma}\right)\right|_{I_n^j}\right]+(2 \sigma)^n c_n\right\}\label{c_n+1}
        \end{align}
    
    \end{subequations}
    \end{proof}

\subsection{Proof of Corollary \ref{Cor c_n+1}}

\begin{proof}
    
From (\ref{fniPos}) and (\ref{f^1_n+1}) we get
 
    \begin{equation}
f_{n+1}^1(u)=\left\{ \frac{e^{\frac{(u-\mu)}{\sigma}}}{2\sigma} c_n +
\sum _{j=1}^{n+1} \frac{e^{\frac{(u-\mu)}{\sigma}}}{(2 \sigma)^{n+1}}\sum_{k=0}^{m^j_n-1}  \int_{I_{n}^j} \left[a_{n}^{(j,k)}(y-n\mu)^k +b_{n}^{(j,k)}(y-n\mu)^ke^{-\frac{2y}{\sigma}} \right] dy\right \} \mathbbm{1}_{I_{n+1}^1}(u)
\end{equation}
    Computing the integrals we obtain 
  \begin{subequations}\label{f1_n+1}
 \begin{align}
f_{n+1}^1(u)&= \left\{\frac{e^{\frac{(u-\mu)}{\sigma}}}{2\sigma} c_n       \right.\\
&+\left.\sum _{j=1}^{n+1} \frac{e^{\frac{(u-\mu)}{\sigma}}}{(2 \sigma)^{n+1}}\sum_{k=0}^{m^j_n-1}   \left\{a_{n}^{(j,k)}\left.\frac{(y-n\mu)^{k+1}}{k+1}\right|_{I_{n}^j}-b_{n}^{(j,k)} \frac{e^{-\frac{2n\mu}{\sigma}}\sigma^{k+1} (-1)^k}{2^{k+1}}\left.\Gamma\left(1+k,\frac{2(y-n\mu)}{\sigma}\right)\right|_{I_{n}^j}\right\} \right \} \mathbbm{1}_{I_{n+1}^1}(u)
\end{align}
\end{subequations}  
 Computing (\ref{f1_n+1}) in $u=0$ and comparing it with (\ref{c_n+1}) we get the thesis. 
\end{proof}

\subsection{Proof of Theorem \ref{Theorem Position -x<mu<0}}
 \begin{proof}
In analogy with the case $\mu>0$, we proceed by induction. 
 
 {\bf Case n=1}
 This part of the proof coincides with the analogous part of the proof of Theorem \ref{Theorem Position mu>=0}. Using (\ref{f_1 generale}) according to $x+\mu$ positive or negative, we get (\ref{fn_b}).

 {\bf Case n} Let us assume that (\ref{fn_b}) holds for $n$, and we show that they  hold for $n+1$, for $n\geq 1$. \\
 In analogy with the proof of Theorem \ref{Theorem Position mu>=0}, using the transition density function (\ref{denstransiz}) in (\ref{CK1}), we get
 \begin{align}
 f_{n+1}(u)&=\mathbbm{1}_{(0,+\infty)}(u)\int_{-\mu}^{+\infty} \frac{e^{-\frac{|u-y-\mu|}{\sigma}}}{2\sigma} f_n(y) dy
 + \delta(u)\int_{-\mu}^{+\infty} \frac{e^{-\frac{y+\mu}{\sigma}}}{2} f_n(y) dy\\
 &+\mathbbm{1}_{(0,+\infty)}(u)\int_0^{-\mu} \frac{e^{-\frac{u-y-\mu}{\sigma}}}{2\sigma} f_n(y) dy
 + \delta(u)\int_0^{-\mu}\left(1-\frac{e^{\frac{y+\mu}{\sigma}}}{2}\right) f_n(y) dy \nonumber
 \end{align}

Since, by induction, $f_n(u)$ is given by (\ref{fn_b}), for $u>0$ we get
\begin{align} \label{partecont}
 f_{n+1}(u)&=\left\{\sum_{i=1}^2\left[\int_0^{-\mu}\frac{e^{-\frac{u-y-\mu}{\sigma}}}{2\sigma} f_n^i(y) dy+\int_{-\mu}^{+\infty} \frac{e^{-\frac{|u-y-\mu|}{\sigma}}}{2\sigma} f_n^i(y) dy \right] + \left[\frac{e^{-\frac{u-\mu}{\sigma}}}{2\sigma} f_n^0(0)\right]\right\}\mathbbm{1}_{(0,+\infty)}(u)
 \end{align}
while for $u=0$ we have
\begin{align}\label{fdirac 3.4}
 f_{n+1}(0)&=  \sum_{i=1}^2\left[\int_0^{-\mu}\left(1-\frac{e^{\frac{y+\mu}{\sigma}}}{2}\right) f_n^i(y) dy+\int_{-\mu}^{+\infty} \frac{e^{-\frac{y+\mu}{\sigma}}}{2} f_n^i(y) dy\right]+\left[\left(1-\frac{e^{\frac{\mu}{\sigma}}}{2}\right) f_n^0(0)\right]
 \end{align}
 Let us focus on the continuous part (\ref{partecont}). Expanding the modulus we get
\begin{align}\label{f pre divisione}
 f_{n+1}(u)=\left\{\sum_{i=1}^2\left[\int_{0}^{u-\mu} \frac{e^{-\frac{u-y-\mu}{\sigma}}}{2\sigma} f_n^i(y) dy+\int_{u-\mu}^{+\infty} \frac{e^{\frac{u-y-\mu}{\sigma}}}{2\sigma} f_n^i(y) dy \right] + \frac{e^{-\frac{u-\mu}{\sigma}}}{2\sigma} f_n^0(0)\right\}\mathbbm{1}_{(0,+\infty)}(u)
\end{align} 

Now, in order to complete the proof we have to distinguish two cases according to the sign of $x+n\mu$.

Let us consider $x+n\mu>0$. Observe that the partition of step $n$ (\ref{partitionTh2}) is updated in step $n+1$ after checking whether $u-\mu<x+n\mu$. The corresponding functions $f_{n+1}^i$, $i=1,2$ are
 \begin{align}\label{f_n+1^1}
 f_{n+1}^1(u)&=\left\{\int_{0}^{u-\mu} \frac{e^{-\frac{u-y-\mu}{\sigma}}}{2\sigma} f_n^1(y) dy+  \int_{u-\mu}^{x+n\mu} \frac{e^{\frac{u-y-\mu}{\sigma}}}{2\sigma} f_n^1(y) dy+\int_{x+n\mu}^{+\infty} \frac{e^{\frac{u-y-\mu}{\sigma}}}{2\sigma} f_n^2(y) dy  + \frac{e^{\frac{-u+\mu}{\sigma}}}{2\sigma} f_n^0(0)\right\}\mathbbm{1}_{I_{n+1}^1}(u)
\end{align} 
 \begin{align}\label{f_n+1^2}
 f_{n+1}^2(u)&=\left\{\int_{0}^{x+n\mu} \frac{e^{-\frac{u-y-\mu}{\sigma}}}{2\sigma} f_n^1(y) dy+  \int_{x+n\mu}^{u-\mu} \frac{e^{\frac{u-y-\mu}{\sigma}}}{2\sigma} f_n^2(y) dy+\int_{u-\mu}^{+\infty} \frac{e^{-\frac{u-y-\mu}{\sigma}}}{2\sigma} f_n^2(y) dy  + \frac{e^{\frac{-u+\mu}{\sigma}}}{2\sigma} f_n^0(0)\right\}\mathbbm{1}_{I_{n+1}^2}(u)
\end{align} 
Substituting in (\ref{f_n+1^1}) the expressions given in (\ref{f_n T2})
 \begin{align}
f_{n+1}^1(u)&= \left\{\frac{e^{-\frac{(u-\mu)}{\sigma}}}{(2 \sigma)^{n+1}}\sum_{k=0}^{n-1} \int_{0}^{u-\mu}  \left[a_{n}^{(1,k)}(y-n\mu)^k e^{\frac{2y}{\sigma}}+b_{n}^{(1,k)}(y-n\mu)^k\right]dy\right.\nonumber\\
&+\frac{e^{\frac{(u-\mu)}{\sigma}}}{(2 \sigma)^{n+1}}\sum_{k=0}^{n-1} \int_{u-\mu}^{x+n\mu} \left[a_{n}^{(1,k)}(y-n\mu)^k +b_{n}^{(1,k)}(y-n\mu)^k e^{-\frac{2y}{\sigma}}\right] dy\nonumber\\
&+\frac{e^{\frac{(u-\mu)}{\sigma}}}{(2 \sigma)^{n+1}}\sum_{k=0}^{n-1} \int_{x+n\mu}^{\infty}\left[a_{n}^{(2,k)}(y-n\mu)^k +b_{n}^{(2,k)}(y-n\mu)^k e^{-\frac{2y}{\sigma}}\right] dy\nonumber \\
&\left.+ \frac{e^{\frac{-u+\mu}{\sigma}}}{2\sigma} f_n^0(0) \right \} \mathbbm{1}_{I_{n+1}^1}(u)
\end{align}

Computing the integrals we get
 \begin{align}\label{f^1 T2}
f_{n+1}^1(u)
&= \left\{\frac{e^{-\frac{(u-\mu)}{\sigma}}}{(2 \sigma)^{n+1}}\sum_{k=0}^{n-1} \left\{a_{n}^{(1,k)}\frac{e^{\frac{2n\mu}{\sigma}}\sigma^{k+1} (-1)^k}{2^{k+1}}\left[\left.\Gamma\left(1+k,\frac{-2(y-n\mu)}{\sigma}\right)\right|_{0}^{u-\mu} \right]+b_{n}^{(1,k)}\left.\frac{(y-n\mu)^{k+1}}{k+1}\right|_{0}^{u-\mu} \right\}\right.\nonumber\\
&+\frac{e^{\frac{(u-\mu)}{\sigma}}}{(2 \sigma)^{n+1}}\sum_{k=0}^{n-1} \left\{a_{n}^{(1,k)}\left.\frac{(y-n\mu)^{k+1}}{k+1}\right|_{u-\mu}^{x+n\mu}-b_{n}^{(1,k)}\frac{e^{\frac{-2n\mu}{\sigma}}\sigma^{k+1}}{2^{k+1}}\left[\left.\Gamma\left(1+k,\frac{2(y-n\mu)}{\sigma}\right)\right|_{u-\mu}^{x+n\mu}\right]
\right\}\nonumber\\
&+\frac{e^{\frac{(u-\mu)}{\sigma}}}{(2 \sigma)^{n+1}}\sum_{k=0}^{n-1} \left\{a_{n}^{(2,k)}\left.\frac{(y-n\mu)^{k+1}}{k+1}\right|_{x+n\mu}^{ \infty}+b_{n}^{(2,k)}\frac{e^{\frac{-2n\mu}{\sigma}}\sigma^{k+1}}{2^{k+1}}\left[\Gamma\left(1+k,\frac{2x}{\sigma}\right)\right]
\right\} \nonumber\\
&\left.+ \frac{e^{\frac{-u+\mu}{\sigma}}}{2\sigma} f_n^0(0) \right \} \mathbbm{1}_{I_{n+1}^1}(u)
\end{align}
Proceeding in an analogous way for (\ref{f_n+1^2}) we compute 

 \begin{align}
f_{n+1}^2(u)&= \left\{\frac{e^{-\frac{(u-\mu)}{\sigma}}}{(2 \sigma)^{n+1}}\sum_{k=0}^{n-1} \int_{0}^{x+n\mu}  \left[a_{n}^{(1,k)}(y-n\mu)^k e^{\frac{2y}{\sigma}}+b_{n}^{(1,k)}(y-n\mu)^k\right]dy\right.\nonumber\\
&+\frac{e^{-\frac{(u-\mu)}{\sigma}}}{(2 \sigma)^{n+1}}\sum_{k=0}^{n-1} \int_{x+n\mu}^{u-\mu} \left[a_{n}^{(2,k)}(y-n\mu)^ke^{\frac{2y}{\sigma}} +b_{n}^{(2,k)}(y-n\mu)^k \right] dy\nonumber\\
&+\frac{e^{\frac{(u-\mu)}{\sigma}}}{(2 \sigma)^{n+1}}\sum_{k=0}^{n-1} \int_{u-\mu}^{\infty}\left[a_{n}^{(2,k)}(y-n\mu)^k +b_{n}^{(2,k)}(y-n\mu)^k e^{-\frac{2y}{\sigma}}\right] dy \left.+ \frac{e^{\frac{-u+\mu}{\sigma}}}{2\sigma} f_n^0(0) \right \} \mathbbm{1}_{I_{n+1}^2}(u)
\end{align}
Computing the integrals we get
 \begin{align}\label{f^2 T2}
f_{n+1}^2(u)
&= \left\{\frac{e^{-\frac{(u-\mu)}{\sigma}}}{(2 \sigma)^{n+1}}\sum_{k=0}^{n-1} \left\{a_{n}^{(1,k)}\frac{e^{\frac{2n\mu}{\sigma}}\sigma^{k+1} (-1)^k}{2^{k+1}}\left[\left.\Gamma\left(1+k,\frac{-2(y-n\mu)}{\sigma}\right)\right|_{0}^{x+n\mu} \right]+b_{n}^{(1,k)}\left.\frac{(y-n\mu)^{k+1}}{k+1}\right|_{0}^{x+n\mu} \right\}\right.\nonumber\\
&+\frac{e^{-\frac{(u-\mu)}{\sigma}}}{(2 \sigma)^{n+1}}\sum_{k=0}^{n-1} \left\{ a_{n}^{(2,k)}\frac{e^{\frac{2n\mu}{\sigma}}\sigma^{k+1} (-1)^k}{2^{k+1}}\left[\left.\Gamma\left(1+k,\frac{-2(y-n\mu)}{\sigma}\right)\right|_{x+n\mu}^{u-\mu}\right] +b_{n}^{(2,k)}\left.\frac{(y-n\mu)^{k+1}}{k+1}\right|_{x+n\mu}^{u-\mu}\right\}\nonumber\\
&+\frac{e^{\frac{(u-\mu)}{\sigma}}}{(2 \sigma)^{n+1}}\sum_{k=0}^{n-1} \left\{a_{n}^{(2,k)}\left.\frac{(y-n\mu)^{k+1}}{k+1}\right|_{u-\mu}^{ \infty}+b_{n}^{(2,k)}\frac{e^{\frac{-2n\mu}{\sigma}}\sigma^{k+1}}{2^{k+1}}\left[\Gamma\left(1+k,\frac{2(u-(n+1)\mu)}{\sigma}\right)\right]
\right\} \nonumber\\
&\left.+ \frac{e^{\frac{-u+\mu}{\sigma}}}{2\sigma} f_n^0(0) \right \} \mathbbm{1}_{I_{n+1}^2}(u)
\end{align}

Expanding the incomplete Gamma function (\ref{gammaExpansion}) in (\ref{f^1 T2}) and (\ref{f^2 T2}) and collecting powers $(u-(n+1)\mu)^je^{\frac{u}{\sigma}}$ and $(u-(n+1)\mu)^je^{-\frac{u}{\sigma}}$, $j=0,\dots, n-1$ we recognize the coefficients (\ref{coeff T2}).

Observe that the divergence of $(y-n\mu)^{k+1}$ on $y=+\infty$ in (\ref{f^1 T2}) and (\ref{f^2 T2}) is not be a problem since, from the recurrent expression (\ref{ricorsione a}) of the coefficients and (\ref{a_1^(2,0)}), $a_2^{(n+1,k)}\equiv 0$ for each admissible $k$.

Let us now consider $c_{n+1}=f_{n+1}(0)$, given by (\ref{fdirac 3.4}). This value changes according to the value of $x$ with respect to the following intervals $[0, -(n+1)\mu)$ and $[ -(n+1)\mu, \infty)$.
Making use of the indicators of such intervals, we get
\begin{subequations} \label{f_0 T3}
\begin{align}
 c_{n+1}&=\left\{ \int_0^{x+n\mu}\left(1-\frac{e^{\frac{y+\mu}{\sigma}}}{2}\right) f_n^1(y) dy+\int_{x+n\mu}^{-\mu}\left(1-\frac{e^{\frac{y+\mu}{\sigma}}}{2}\right) f_n^2(y) dy+\int_{-\mu}^{+\infty} \frac{e^{-\frac{y+\mu}{\sigma}}}{2} f_n^2(y) dy\right.\nonumber  \label{f_0 left} \\ 
 &+  \left.\left(1-\frac{e^{\frac{\mu}{\sigma}}}{2}\right) f_n^0(0)\right\}\mathbbm{1}_{[0,-(n+1)\mu)}(x) \\
&+ \left\{\int_0^{-\mu}\left(1-\frac{e^{\frac{y+\mu}{\sigma}}}{2}\right) f_n^1(y) dy+\int_{-\mu}^{x+n\mu}\frac{e^{-\frac{y+\mu}{\sigma}}}{2} f_n^1(y) dy+\int_{x+n\mu}^{+\infty} \frac{e^{-\frac{y+\mu}{\sigma}}}{2} f_n^2(y) dy\right.\nonumber\\\label{f_0 right}
 &+  \left.\left(1-\frac{e^{\frac{\mu}{\sigma}}}{2}\right) f_n^0(0)\right\}  \mathbbm{1}_{[-(n+1)\mu),\infty)}(x)\\
 &=\tilde{f}_{n+1}(0)+\tilde{\tilde{f}}_{n+1}(0)\nonumber
 \end{align}
\end{subequations}
where  $\tilde{f}_{n+1}(0)$ corresponds to (\ref{f_0 left}) and $\tilde{\tilde{f}}_{n+1}(0)$ corresponds to (\ref{f_0 right}).

Substituting (\ref{f_n T2}) in (\ref{f_0 T3}) we get

\begin{align}
 \tilde{f}_{n+1}(0)&=\left\{\frac{1}{(2\sigma)^n}\sum_{k=0}^{n-1}\int_0^{x+n\mu}\left[a_{n}^{(1,k)}(y-n\mu)^k e^{\frac{y}{\sigma}}+b_{n}^{(1,k)}(y-n\mu)^k e^{-\frac{y}{\sigma}}\right]dy\right.\\
&-\frac{e^{\frac{\mu}{\sigma}}}{2^{n+1}\sigma^n}\sum_{k=0}^{n-1}\int_0^{x+n\mu}\left[a_{n}^{(1,k)}(y-n\mu)^k e^{\frac{2y}{\sigma}}+b_{n}^{(1,k)}(y-n\mu)^k\right]dy\nonumber\\
 &+\frac{1}{(2\sigma)^n}\sum_{k=0}^{n-1}\int_{x+n\mu}^{-\mu}\left[a_{n}^{(2,k)}(y-n\mu)^k e^{\frac{y}{\sigma}}+b_{n}^{(2,k)}(y-n\mu)^k e^{-\frac{y}{\sigma}}\right] dy\nonumber\\
&-\frac{e^{\frac{\mu}{\sigma}}}{2^{n+1}\sigma^n}\sum_{k=0}^{n-1}\int_{x+n\mu}^{-\mu}\left[a_{n}^{(2,k)}(y-n\mu)^k e^{\frac{2y}{\sigma}}+b_{n}^{(2,k)}(y-n\mu)^k\right] dy\nonumber\\
&+\frac{e^{-\frac{\mu}{\sigma}}}{2^{n+1}\sigma^n}\sum_{k=0}^{n-1}\int_{-\mu}^{+\infty} \left[a_{n}^{(2,k)}(y-n\mu)^k +b_{n}^{(2,k)}(y-n\mu)^k e^{-\frac{2y}{\sigma}}\right]dy+  \left.\left(1-\frac{e^{\frac{\mu}{\sigma}}}{2}\right) c_n\right\} \mathbbm{1}_{[0,-(n+1)\mu)}(x)\nonumber
 \end{align}

\begin{align}
 \tilde{\tilde{f}}_{n+1}(0)&=\left\{\frac{1}{(2\sigma)^n}\sum_{k=0}^{n-1}\int_0^{-\mu}\left[a_{n}^{(1,k)}(y-n\mu)^k e^{\frac{y}{\sigma}}+b_{n}^{(1,k)}(y-n\mu)^k e^{-\frac{y}{\sigma}}\right]dy\right.\\
 &-\frac{e^{\frac{\mu}{\sigma}}}{2^{n+1}\sigma^n}\sum_{k=0}^{n-1}\int_0^{-\mu}\left[a_{n}^{(1,k)}(y-n\mu)^k e^{\frac{2y}{\sigma}}+b_{n}^{(1,k)}(y-n\mu)^k \right]dy\nonumber\\
 &+\frac{e^{-\frac{\mu}{\sigma}}}{2^{n+1}\sigma^n}\sum_{k=0}^{n-1}\int_{-\mu}^{x+n\mu}\left[a_{n}^{(1,k)}(y-n\mu)^k+b_{n}^{(1,k)}(y-n\mu)^k e^{-\frac{2y}{\sigma}}\right] dy\nonumber\\
&+\frac{e^{-\frac{\mu}{\sigma}}}{2^{n+1}\sigma^n}\sum_{k=0}^{n-1}\int_{x+n\mu}^{+\infty}  \left[a_{n}^{(2,k)}(y-n\mu)^k+b_{n}^{(2,k)}(y-n\mu)^k e^{-\frac{2y}{\sigma}}\right]dy+ \left. \left[\left(1-\frac{e^{\frac{\mu}{\sigma}}}{2}\right) c_n\right] \right\}\mathbbm{1}_{[-(n+1)\mu),\infty)}(x)\nonumber
 \end{align}

Computing the integrals we get the thesis (\ref{massa pos}) and 
(\ref{massa neg}).

Let us consider now $x+n\mu\leq0$. 

Note that in this case $I_n^1= \emptyset$ and $f_n^1$ is identically $0$ while $f_{n+1}^2(u)$ becomes
 \begin{align}\label{f_n+1^2 bis}
 f_{n+1}^2(u)&=\left\{ \int_0^{u-\mu} \frac{e^{\frac{u-y-\mu}{\sigma}}}{2\sigma} f_n^2(y) dy+\int_{u-\mu}^{+\infty} \frac{e^{-\frac{u-y-\mu}{\sigma}}}{2\sigma} f_n^2(y) dy  + \frac{e^{\frac{-u+\mu}{\sigma}}}{2\sigma} f_n^0(0)\right\}\mathbbm{1}_{I_{n+1}^2}(u)
\end{align}

Substituting in (\ref{f_n+1^2 bis}) the expressions given in (\ref{fni 3.4})

 \begin{align}
f_{n+1}^2(u)&= \left\{\frac{1}{(2 \sigma)^n}\sum_{k=0}^{n-1} \int_{0}^{u-\mu}
\frac{e^{-\frac{(u-y-\mu)}{\sigma}}}{2\sigma} \left[a_{n}^{(2,k)}(y-n\mu)^k e^{\frac{y}{\sigma}}+b_{n}^{(2,k)}(y-n\mu)^k e^{-\frac{y}{\sigma}}\right]\right. dy \nonumber\\
&+\frac{1}{(2 \sigma)^n}\sum_{k=0}^{n-1} \int_{u-\mu}^{\infty}\frac{e^{\frac{(u-y-\mu)}{\sigma}}}{2\sigma} \left[a_{n}^{(2,k)}(y-n\mu)^k e^{\frac{y}{\sigma}}+b_{n}^{(2,k)}(y-n\mu)^k e^{-\frac{y}{\sigma}}\right] dy \nonumber\\
&\left.+ \frac{e^{\frac{-u+\mu}{\sigma}}}{2\sigma} f_n^0(0) \right \} \mathbbm{1}_{I_{n+1}^2}(u)
\end{align}
Computing the integrals we get
 \begin{align}
f_{n+1}^2(u)
&= \left\{\frac{e^{-\frac{(u-\mu)}{\sigma}}}{(2 \sigma)^{n+1}}\sum_{k=0}^{n-1} \left\{ a_{n}^{(2,k)}\frac{e^{\frac{2n\mu}{\sigma}}\sigma^{k+1} (-1)^k}{2^{k+1}}\left[\left.\Gamma\left(1+k,\frac{-2(y-n\mu)}{\sigma}\right)\right|_{0}^{u-\mu}\right] +b_{n}^{(2,k)}\left.\frac{(y-n\mu)^{k+1}}{k+1}\right|_{0}^{u-\mu}\right\} \right.\nonumber\\
&+\frac{e^{\frac{(u-\mu)}{\sigma}}}{(2 \sigma)^{n+1}}\sum_{k=0}^{n-1} \left\{a_{n}^{(2,k)}\left.\frac{(y-n\mu)^{k+1}}{k+1}\right|_{u-\mu}^{ \infty}+b_{n}^{(2,k)}\frac{e^{\frac{-2n\mu}{\sigma}}\sigma^{k+1}}{2^{k+1}}\left[\Gamma\left(1+k,\frac{2(u-(n+1)\mu)}{\sigma}\right)\right]
\right\} \nonumber\\
&\left.+ \frac{e^{\frac{-u+\mu}{\sigma}}}{2\sigma} f_n^0(0) \right \} \mathbbm{1}_{I_{n+1}^2}(u)
\end{align}

The coefficients of $e^{\frac{u}{\sigma}}$ and  $e^{-\frac{u}{\sigma}}$ 
give (\ref{coeff a b T2bis}).

It remains to compute the mass in $0$ when $x+n\mu\leq0$. Note that in this setting we have only to consider the case $x+n\mu\leq -\mu$ since $x+n\mu> -\mu$ is impossible.\\
For $x<-(n+1)\mu$ (\ref{fdirac 3.4}) becomes

\begin{align}
 c_{n+1}=f_{n+1}(0)&= \left[\int_{0}^{-\mu}\left(1-\frac{e^{\frac{y+\mu}{\sigma}}}{2}\right) f_n^2(y) dy+\int_{-\mu}^{+\infty} \frac{e^{-\frac{y+\mu}{\sigma}}}{2} f_n^2(y) dy\right]   +  \left[\left(1-\frac{e^{\frac{\mu}{\sigma}}}{2}\right) f_n^0(0)\right]
 \end{align}

Substituting (\ref{fni 3.4}) 
\begin{subequations}
\begin{align}
c_{n+1} &=\frac{1}{(2\sigma)^n}\sum_{k=0}^{n-1}\int_{0}^{-\mu}\left[a_{n}^{(2,k)}(y-n\mu)^k e^{\frac{y}{\sigma}}+b_{n}^{(2,k)}(y-n\mu)^k e^{-\frac{y}{\sigma}}\right] dy\\
&-\frac{e^{\frac{\mu}{\sigma}}}{2^{n+1}\sigma^n}\sum_{k=0}^{n-1}\int_{0}^{-\mu}\left[a_{n}^{(2,k)}(y-n\mu)^k e^{\frac{2y}{\sigma}}+b_{n}^{(2,k)}(y-n\mu)^k\right] dy\\
&+\frac{e^{-\frac{\mu}{\sigma}}}{2^{n+1}\sigma^n}\sum_{k=0}^{n-1}\int_{-\mu}^{+\infty} \left[a_{n}^{(2,k)}(y-n\mu)^k +b_{n}^{(2,k)}(y-n\mu)^k e^{-\frac{2y}{\sigma}}\right]dy\\
 &+  \left[\left(1-\frac{e^{\frac{\mu}{\sigma}}}{2}\right) c_n\right]
 \end{align}
\end{subequations}
that gives the result (\ref{c_n+1 bis}).
 \end{proof}

\subsection{Proof of Theorem \ref{Theorem Position mu<=-x}}
    \begin{proof}

 {\bf Case n=1}
 
  This part of the proof coincides with the analogous part of the proof of Theorem \ref{Theorem Position mu>=0}.

 {\bf Case n}

Let us assume that (\ref{fn_1 3.5}) holds for $n$, and we show that it holds for $n+1$, for $n\geq 1$. 

 In analogy with the proof of Theorem \ref{Theorem Position mu>=0}, using the transition density function (\ref{denstransiz}) in (\ref{CK1}), we get

 \begin{align}
 f_{n+1}(u)&=\mathbbm{1}_{(0,+\infty)}(u)\int_{-\mu}^{+\infty} \frac{e^{-\frac{|u-y-\mu|}{\sigma}}}{2\sigma} f_n(y) dy
 + \delta(u)\int_{-\mu}^{+\infty} \frac{e^{-\frac{y+\mu}{\sigma}}}{2} f_n(y) dy\\
 &+\mathbbm{1}_{(0,+\infty)}(u)\int_0^{-\mu} \frac{e^{-\frac{u-y-\mu}{\sigma}}}{2\sigma} f_n(y) dy
 + \delta(u)\int_0^{-\mu}\left(1-\frac{e^{\frac{y+\mu}{\sigma}}}{2}\right) f_n(y) dy \nonumber
 \end{align}

Using (\ref{fn 3.5}) we get
 \begin{align}
 f_{n+1}(u)&=\sum _{i=0}^{1} \left[\int_{-\mu}^{+\infty} \frac{e^{-\frac{|u-y-\mu|}{\sigma}}}{2\sigma} f_n^i(y) dy\right]\mathbbm{1}_{(0,+\infty)}(u)+ \sum _{i=0}^{1} \left[\int_{-\mu}^{+\infty} \frac{e^{-\frac{y+\mu}{\sigma}}}{2} f_n^i(y) dy\right]\delta(u)
\\&\qquad +\sum _{i=0}^{1} \left[\int_0^{-\mu}\frac{e^{-\frac{u-y-\mu}{\sigma}}}{2\sigma} f_n^i(y) dy\right]\mathbbm{1}_{(0,+\infty)}(u)+ \sum _{i=0}^{1} \left[\int_0^{-\mu}\left(1-\frac{e^{\frac{y+\mu}{\sigma}}}{2}\right) f_n^i(y) dy\right]\delta(u)\nonumber
 \end{align}

When $u>0$ we get
\begin{align}
 f_{n+1}(u)&=\left\{\int_{-\mu}^{+\infty} \frac{e^{-\frac{|u-y-\mu|}{\sigma}}}{2\sigma} f_n^1(y) dy+ \int_0^{-\mu}\frac{e^{-\frac{u-y-\mu}{\sigma}}}{2\sigma} f_n^1(y) dy + \frac{e^{-\frac{u-\mu}{\sigma}}}{2\sigma} f_n^0(0)\right\}\mathbbm{1}_{(0,+\infty)}(u)\\
 &=\left\{\int_{0}^{u-\mu} \frac{e^{-\frac{u-y-\mu}{\sigma}}}{2\sigma} f_n^1(y) dy+\int_{u-\mu}^{+\infty} \frac{e^{\frac{u-y-\mu}{\sigma}}}{2\sigma} f_n^1(y) dy + \frac{e^{\frac{-u+\mu}{\sigma}}}{2\sigma} f_n^0(0)\right\}\mathbbm{1}_{(0,+\infty)}(u)\nonumber
 \end{align}

Applying the inductive hypothesis (\ref{f_n T3})
 
\begin{align}
 f_{n+1}(u)
 &=\left\{\frac{1}{(2\sigma)^n}\frac{e^{-\frac{u-\mu}{\sigma}}}{2\sigma} \sum_{j=0}^{n-1}\left(\int_{0}^{u-\mu} b_n^{(1,j)}(y-(n-1)\mu)^jdy\right)  
 \right.\nonumber\\
 &+\frac{1}{(2\sigma)^n}\frac{e^{\frac{u-\mu}{\sigma}}}{2\sigma}\sum_{j=0}^{n-1}\left( \int_{u-\mu}^{+\infty}b_n^{(1,j)}(y-(n-1)\mu)^j e^{-\frac{2y}{\sigma}}dy\right) + \left. \frac{e^{\frac{-u+\mu}{\sigma}}}{2\sigma} f_n^0(0)\right\}\mathbbm{1}_{(0,+\infty)}(u)
\end{align} 

Computing the integrals and expanding the incomplete Gamma functions (\ref{gammaExpansion}) we obtain
\begin{align}
 f_{n+1}(u)
&=\left\{\frac{1}{(2\sigma)^n}\frac{e^{-\frac{u-\mu}{\sigma}}}{2\sigma} \sum_{k=0}^{n-1}\left[\frac{b_n^{(1,k)}(u-n\mu)^{k+1}}{k+1}-\frac{b_n^{(1,k)}(-(n-1)\mu)^{k+1}}{k+1}\right]\right.\nonumber\\
 &+\frac{1}{(2\sigma)^n}\frac{e^{\frac{u-\mu}{\sigma}}}{2\sigma}\sum_{k=0}^{n-1}b_n^{(1,k)}\frac{e^{-\frac{2(n-1)\mu}{\sigma}}\sigma^{k+1}}{2^{k+1}}\left[k! e^{-\frac{2(u-n\mu)}{\sigma}} \sum_{r=0}^k \frac{\left(\frac{2(u-n\mu)}{\sigma}\right)^r}{r!}\right]+ \left. \frac{e^{\frac{-u+\mu}{\sigma}}}{2\sigma} f_n^0(0)\right\}\mathbbm{1}_{(0,+\infty)}(u)
\end{align} 
where we can recognize the recurrence relations (\ref{abc T3})

Let us now work out the calculation of the probability to be in $0$ at step $n+1$

\begin{align}\label{fdirac_c}
 c_{n+1}=f_{n+1}(0)&=  \left[\int_{-\mu}^{+\infty} \frac{e^{-\frac{y+\mu}{\sigma}}}{2} f_n^1(y) dy\right]+  \left[\int_0^{-\mu}\left(1-\frac{e^{\frac{y+\mu}{\sigma}}}{2}\right) f_n^1(y) dy\right] +  \left[\left(1-\frac{e^{\frac{\mu}{\sigma}}}{2}\right) f_n^0(0)\right]
 \end{align}
Substituting in (\ref{fdirac_c}) the expression of the inductive hypothesis (\ref{fn_1 3.5}) we get
\begin{subequations}
\begin{align}
 c_{n+1}
  &=  \left[\frac{e^{-\frac{\mu}{\sigma}}}{2}\int_{-\mu}^{+\infty} \left[ \sum_{k=0}^{n-1}\frac{1}{(2\sigma)^n}\left( b_n^{(1,k)}(y-(n-1)\mu)^ke^{-\frac{2y}{\sigma}} \right) \right] dy\right] \\
  &+  \left[\int_{0}^{-\mu}\sum_{k=0}^{n-1}\frac{1}{(2\sigma)^n}\left( b_n^{(1,k)}(y-(n-1)\mu)^k e^{-\frac{y}{\sigma}} \right)dy\right]-\frac{e^{\frac{\mu}{\sigma}}}{2} \left[\int_{0}^{-\mu}\sum_{k=0}^{n-1}\frac{1}{(2\sigma)^n}\left( b_n^{(1,k)}(y-(n-1)\mu)^k \right) dy\right]\\
  &+  \left[\left(1-\frac{e^{\frac{\mu}{\sigma}}}{2}\right) f_n^0(0)\right]\\
  &=-\frac{e^{-\frac{\mu}{\sigma}}}{2^{n+1}\sigma^{n}}\sum_{k=0}^{n-1}b_n^{(1,k)}\frac{e^{-\frac{2(n-1)\mu}{\sigma}}\sigma^{k+1}}{2^{k+1}}\left[\left.\Gamma\left(1+k,\frac{2(y-(n-1)\mu)}{\sigma}\right)\right|^\infty_{-\mu}\right]\\
  &-\frac{1}{(2\sigma)^{n}}\sum_{k=0}^{n-1}b_n^{(1,k)}e^{-\frac{(n-1)\mu}{\sigma}}\sigma^{k+1}\left[\left.\Gamma\left(1+k,\frac{(y-(n-1)\mu)}{\sigma}\right)\right|_0^{-\mu}\right]\\
  &-\frac{e^{\frac{\mu}{\sigma}}}{2^{n+1}\sigma^{n}}\sum_{k=0}^{n-1}b_n^{(1,k)}\left[\left.\frac{(y-(n-1)\mu)^{k+1}}{k+1}\right|_0^{-\mu}\right]\\
  &+  \left[\left(1-\frac{e^{\frac{\mu}{\sigma}}}{2}\right) f_n^0(0)\right]
 \end{align}
 \end{subequations}
 that gives the coefficient (\ref{c_n+1 T3}).

 \end{proof}   

\section{Proofs of Theorems on the FET}

\subsection{Proof of Theorem \ref{Theorem FPT mu>0 h>mu}}
\begin{proof}
We proceed by induction.

{\bf Case n=1}
Since $P(1|x)=1-F_1(h|x)$ and $0<\mu<h$ by hypothesis, using (\ref{F1(u|x)}) in Lemma \ref{Lemma1}, for $x\in[0,h)$, we have  
\begin{subequations}
\begin{align} 
P(1|x) &=\frac{1}{2}e^{\frac{\mu-h+x}{\sigma}}\mathbbm{1}_{I_1^1}(x) +     \left(	1-\frac{1}{2}e^{-\frac{\mu-h+x}{\sigma}}\right)  \mathbbm{1}_{I_1^2}(x) \\
&=\sum _{i=1}^{2} P_i(1|x) 
\end{align} 
\end{subequations}
where the partition (\ref{PartitionN}) is 
\begin{align*}
	&I_1^1=[0,h-\mu)\\
	&I_1^2=[h-\mu,h)
\end{align*} 
and we easily recognize the coefficients  (\ref{Ci T FET 1}).

{\bf Step n}

 We assume that (\ref{Pn}) holds for $n$, and we show that it holds for $n+1$, for $n\geq 1$. 

Conditioning on the position reached at time $1$ we get
\begin{equation}\label{smoluk}
 P(n+1|x)=\int_0^{h} P(n|y) \mathbb{P}(W_1\in dy)  \qquad n=1,2, \dots
\end{equation}

Substituting (\ref{f_1 generale}) and using the inductive hypothesis (\ref{Pni}), we obtain:

\begin{align}\label{recursion}
	P(n+1|x)&=P_1(n|0)\frac{e^{-\frac{\mu+x}{\sigma}}}{2}+\sum_{j=1}^{i^*-1}\int_{I^j_n} P_j(n|y)\frac{e^{\frac{y-\mu-x}{\sigma}}}{2\sigma}dy
	+\int_{h-(\ell_n-i^*+1)\mu}^{\min(h,\mu+x)}  P_{i^*}(n|y)\frac{e^{\frac{y-\mu-x}{\sigma}}}{2\sigma}dy\\	
	&+\int_{\min(h,\mu+x)}^{\min(h,h-(\ell_n-i^*)\mu)} P_{i^*}(n|y)\frac{e^{\frac{-y+\mu+x}{\sigma}}}{2\sigma}dy
	+\sum_{j=i^*+1}^{\ell_n}\int_{I^j_n} P_j(n|y) \frac{e^{\frac{-y+\mu+x}{\sigma}}}{2\sigma}dy \nonumber
\end{align}
where $I^{i^*}_n$ is the interval containing $x+\mu$, for $n>1$. Note that the value of $i^*$ depends on $x$, $\mu$ and $n$ and $1\leq i^* \leq \ell_n$.

Now we rewrite $P(n+1|x)$ as a sum of  $\ell_{n+1}$ terms $P_i(n+1|x)$, $i=1,\dots, \ell_{n+1}$ where
\begin{align}\label{P_n+1^i_prima1}
  	P_i(n+1|x)&= \left\{
  	P_1(n|0)\frac{e^{-\frac{\mu+x}{\sigma}}}{2}+\sum_{j=1}^{i^*-1}\int_{I^j_n} P_j(n|y)\frac{e^{\frac{y-\mu-x}{\sigma}}}{2\sigma}dy
	+\int_{h-(\ell_n-i^*+1)\mu}^{\min(h,\mu+x)}  P_{i^*}(n|y)\frac{e^{\frac{y-\mu-x}{\sigma}}}{2\sigma}dy \right.\\	
	&+\left. \int_{\min(h,\mu+x)}^{\min(h,h-(\ell_n-i^*)\mu)} P_{i^*}(n|y)\frac{e^{\frac{-y+\mu+x}{\sigma}}}{2\sigma}dy
	+\sum_{j=i^*+1}^{\ell_n}\int_{I^j_n} P_j(n|y)\frac{e^{\frac{-y+\mu+x}{\sigma}}}{2\sigma}dy  \right \} \mathbbm{1}_{I^{i}_{n+1}}(x). \nonumber
\end{align}
Observe that, when $i=\ell_{n+1}$ we have $h-\mu<x<h$, hence $i^*=\ell_{n}+1$ and (\ref{P_n+1^i_prima1}) reduces to the sum of the first two terms.\\
Now, using the inductive hypothesis (\ref{Pni}), we get

	\begin{align}\label{P_i pregamma}
		P_i(n+1|x)&= \left\{\frac{e^{-\frac{\mu+x}{\sigma}}}{2}
		\left[\sum_{k=0}^{m^1_n-1} \frac{1}{2^n \sigma^{n-1}}\left[\alpha_{n}^{(1,k)}(n\mu)^k +\beta_{n}^{(1,k)}(n\mu)^k\right] + \eta_n^{(1,0)}\right]
		\right.\\
		&+\sum _{j=1}^{i^*-1} \frac{1}{2^n \sigma^{n-1}} \int_{I_n^{j}} \frac{e^{\frac{(y-\mu-x)}{\sigma}}}{2\sigma}\left\{ \sum_{k=0}^{m^j_n-1} \left[\alpha_{n}^{(j,k)}(y+n\mu)^k e^{\frac{y}{\sigma}}+\beta_{n}^{(j,k)}(y+n\mu)^k e^{-\frac{y}{\sigma}}\right]+2^n \sigma^{n-1}\eta_{n}^{(j,0)}\right\} dy\nonumber\\
		&+\frac{1}{2^n \sigma^{n-1}} \int_{h-(\ell_n-i^*+1)\mu}^{\min(h,x+\mu)} \frac{e^{\frac{(y-\mu-x)}{\sigma}}}{2\sigma} \left\{\sum_{k=0}^{m^{i^*}_n-1} \left[\alpha_{n}^{(i^*,k)}(y+n\mu)^k e^{\frac{y}{\sigma}}+\beta_{n}^{(i^*,k)}(y+n\mu)^k e^{-\frac{y}{\sigma}}\right] +\eta_{n}^{(i^*,0)}2^n \sigma^{n-1}\right\}dy\nonumber
   \end{align}
    \begin{align}
		&+\frac{1}{2^n \sigma^{n-1}} \int_{\min(h,x+\mu)}^{\min(h,h-(\ell_n-i^*)\mu)} \frac{e^{\frac{-(y-\mu-x)}{\sigma}}}{2\sigma} \left\{\sum_{k=0}^{m^{i^*}_n-1}\left[\alpha_{n}^{(i^*,k)}(y+n\mu)^k e^{\frac{y}{\sigma}}+\beta_{n}^{(i^*,k)}(y+n\mu)^k e^{-\frac{y}{\sigma}}\right]+\eta_{n}^{(i^*,0)}2^n \sigma^{n-1}\right\} dy\nonumber\\
		&\left.+\sum _{j=i^*+1}^{\ell_n} \frac{1}{2^n \sigma^{n-1}} \int_{I_j} \frac{e^{\frac{-(y-\mu-x)}{\sigma}}}{2\sigma}\left\{\sum_{k=0}^{m^j_n-1} \left[\alpha_{n}^{(j,k)}(y+n\mu)^k e^{\frac{y}{\sigma}}+\beta_{n}^{(j,k)}(y+n\mu)^k e^{-\frac{y}{\sigma}}\right]+\eta_{n}^{(j,0)}2^n \sigma^{n-1}\right\} dy \right \} \mathbbm{1}_{I_{n+1}^i}(x)\nonumber
	\end{align}

Recognizing Gamma functions we get

\begin{align}
	P_i(n+1|x)&=\left\{
	\frac{e^{-\frac{(x+\mu)}{\sigma}}}{2^{n+1} \sigma^n} \sum _{j=0}^{i^*-1} K_n^j\right.+\frac{e^{-\frac{(x+\mu)}{\sigma}}}{2^{n+1} \sigma^n}\left\{\sum_{k=0}^{m^{i^*}_n-1} \left[\alpha_{n}^{(i^*,k)}\frac{e^{-\frac{2n\mu}{\sigma}}\sigma^{k+1} (-1)^k}{2^{k+1}}\left.\Gamma\left(1+k,\frac{-2(y+n\mu)}{\sigma}\right)\right|_{h-(\ell_n-i^*+1)\mu}^{\min(h,x+\mu)} \right.\right.
 \nonumber \\
&\left.+\beta_{n}^{(i^*,k)}\left.\frac{(y+n\mu)^{k+1}}{k+1}\right|_{h-(\ell_n-i^*+1)\mu}^{\min(h,x+\mu)}\right]
	\left.\left. +\eta_{n}^{(i^*,0)}(2\sigma)^n e^{\frac{y}{\sigma}}\right|_{h-(\ell_n-i^*+1)\mu}^{\min(h,x+\mu)} \right\}  \nonumber
   \end{align}
 \begin{align}
	&+\frac{e^{\frac{(x+\mu)}{\sigma}}}{2^{n+1} \sigma^n} \left\{\sum_{k=0}^{m^{i^*}_n-1} \left[\alpha_{n}^{(i^*,k)}\left.\frac{(y+n\mu)^{k+1}}{k+1}\right|_{\min(h,\mu+x)}^{\min(h,h-(\ell_n-i^*)\mu)}\right.\right.\label{termineNotoDx}\nonumber\\
	&-\beta_{n}^{(i^*,k)}\frac{e^{\frac{2n\mu}{\sigma}}\sigma^{k+1} }{2^{k+1}}\left.\left.\Gamma\left(1+k,\frac{2(y+n\mu)}{\sigma}\right)\right|_{\min(h,x+\mu)}^{\min(h,h-(\ell_n-i^*)\mu)} \right]
	\left.\left. -\eta_{n}^{(i^*,0)}(2\sigma)^n e^{-\frac{y}{\sigma}}\right|^{\min(h,h-(\ell_n-i^*)\mu)}_{\min(h,x+\mu)} \right\} \nonumber \\ 
	&\left.+\frac{e^{\frac{(x+\mu)}{\sigma}}}{2^{n+1} \sigma^{n}} \sum _{j=i^*+1}^{\ell_n} K_n^j \right \} \mathbbm{1}_{I_{n+1}^i}(x)
\end{align}
where $K_n^j$ are given by (\ref{K tot}).
Now we can use the expansion (\ref{gammaExpansion}) to get

	\begin{align}
	P_i(n+1|x)&=\left\{
	\frac{e^{-\frac{(x+\mu)}{\sigma}}}{2^{n+1} \sigma^n} \sum _{j=0}^{i^*-1} K_n^j\right.+\frac{e^{-\frac{(x+\mu)}{\sigma}}}{2^{n+1} \sigma^n}\left\{\sum_{k=0}^{m^{i^*}_n-1} \left[\alpha_{n}^{(i^*,k)}\frac{\sigma^{k+1} (-1)^k}{2^{k+1}}\left.k! e^{\frac{2y}{\sigma}} \sum_{r=0}^k \frac{(-2)^r}{r! \sigma^r}(y+n\mu)^r\right|_{h-(\ell_n-i^*+1)\mu}^{\min(h,x+\mu)}\right.\right.\nonumber\\
&\qquad+\left.\beta_{n}^{(i^*,k)}\left.\frac{(y+n\mu)^{k+1}}{k+1}\right|_{h-(\ell_n-i^*+1)\mu}^{\min(h,x+\mu)}
	\right] \left.\left.+\eta_{n}^{(i^*,0)}(2\sigma)^n e^{\frac{y}{\sigma}}\right|_{h-(\ell_n-i^*+1)\mu}^{\min(h,x+\mu)} \right\} \nonumber\\
	&+\frac{e^{\frac{(x+\mu)}{\sigma}}}{2^{n+1} \sigma^n}\left\{\sum_{k=0}^{m^{i^*}_n-1} \left[\alpha_{n}^{(i^*,k)}\frac{(y+n\mu)^{k+1}}{k+1}
	-\beta_{n}^{(i^*,k)}\frac{\sigma^{k+1} }{2^{k+1}}\left.k! e^{\frac{-2y}{\sigma}} \sum_{r=0}^k \frac{(2)^r}{r! \sigma^r}(y+n\mu)^r\right]\right|^{\min(h,h-(\ell_n-i^*)\mu)}_{\min(h,x+\mu)}\right.\nonumber\\
	 &\qquad\left.\left. -\eta_{n}^{(i^*,0)}(2\sigma)^n e^{-\frac{y}{\sigma}}\right|^{\min(h,h-(\ell_n-i^*)\mu)}_{\min(h,x+\mu)} \right\} +\left.\frac{e^{\frac{(x+\mu)}{\sigma}}}{2^{n+1} \sigma^n} \sum _{j=i^*+1}^{\ell_n} K_n^j\right \} \mathbbm{1}_{I_{n+1}^i}(x)
	\end{align}

Substituting the integration values, for $i<\ell_{n+1}$ we get
	\begin{align}
	P_i(n+1|x)&=\left\{
	\frac{e^{-\frac{x}{\sigma}}}{2^{n+1} \sigma^n} e^{-\frac{\mu}{\sigma}}\sum _{j=0}^{i^*-1} K_n^j\right.\label{N1 primo termine}+\frac{e^{-\frac{x}{\sigma}}}{2^{n+1} \sigma^n}e^{-\frac{\mu}{\sigma}}\left\{\sum_{k=0}^{m^{i^*}_n-1} \left\{\alpha_{n}^{(i^*,k)}\frac{\sigma^{k+1} (-1)^k}{2^{k+1}}k!\right.\right.\nonumber\\
	& \qquad\left[ e^{\frac{2(x+\mu)}{\sigma}} \sum_{r=0}^k \frac{(-2)^r}{r! \sigma^r}((x+\mu)+n\mu)^r-e^{\frac{2(h-(\ell_n-i^*+1)\mu)}{\sigma}} \sum_{r=0}^k \frac{(-2)^r}{r! \sigma^r}((h-(\ell_n-i^*+1)\mu)+n\mu)^r \right]\nonumber\\
		&\left.+\beta_{n}^{(i^*,k)}\left[\frac{((x+\mu)+n\mu)^{k+1}}{k+1}-\frac{((h-(\ell_n-i^*+1)\mu)+n\mu)^{k+1}}{k+1}\right]\right\}\nonumber\\
	&+\left.\eta_{n}^{(i^*,0)}(2\sigma)^n\left[e^{\frac{x+\mu}{\sigma}}- e^{\frac{h-(\ell_n-i^*+1)\mu}{\sigma}}\right]\right\} \nonumber\\
	&+\frac{e^{\frac{x}{\sigma}}}{2^{n+1} \sigma^n}e^{\frac{\mu}{\sigma}}\left\{\sum_{k=0}^{m^{i^*}_n-1} \left\{\alpha_{n}^{(i^*,k)}\left[\frac{((h-(\ell_n-i^*)\mu)+n\mu)^{k+1}}{k+1}-\frac{(\mu+x+n\mu)^{k+1}}{k+1}\right]\right.\right.\nonumber\\
	&\left.-\beta_{n}^{(i^*,k)}\frac{\sigma^{k+1} }{2^{k+1}}k!\left[ e^{\frac{-2(h-(\ell_n-i^*)\mu)}{\sigma}} \sum_{r=0}^k \frac{(2)^r}{r! \sigma^r}(h-(\ell_n-i^*)\mu+n\mu)^r- e^{\frac{-2(x+\mu)}{\sigma}} \sum_{r=0}^k \frac{(2)^r}{r! \sigma^r}(x+\mu+n\mu)^r\right]\right\}\nonumber\\
	 &-\left.\eta_{n}^{(i^*,0)}(2\sigma)^n\left[ e^{-\frac{h-(\ell_n-i^*)\mu}{\sigma}}-e^{-\frac{x+\mu}{\sigma}}\right]\right\} +\left.\frac{e^{\frac{x}{\sigma}}}{2^{n+1} \sigma^n}e^{\frac{\mu}{\sigma}} \sum _{j=i^*+1}^{\ell_n} K_n^j \right \} \mathbbm{1}_{I_{n+1}^i}(x)
	\end{align}
 and if $i=\ell_{n+1}$, $P_i(n+1|x)$ reduces to  
\begin{align}
	P_{\ell_{n+1}}(n+1|x)=
	\frac{e^{-\frac{x}{\sigma}}}{2^{n+1} \sigma^n} 
e^{-\frac{\mu}{\sigma}}\sum _{j=0}^{\ell_n} K_n^j.
 \end{align}
Grouping out the terms $e^{\frac{x}{\sigma}}(x+(n+1)\mu)^j$ and  $e^{-\frac{x}{\sigma}}(x+(n+1)\mu)^j$, we get the coefficients $\alpha_{n+1}^{(i,j)}$ and $\beta_{n+1}^{(i,j)}$ respectively.
Furthermore (\ref{eta}) is obtained grouping constant terms.
\end{proof}

\subsection{Proof of Corollary \ref{Theorem FPT mu>0 0<h<mu}}
\begin{proof}
We proceed by induction. However, in this case we need two steps before starting the induction. 
 
 {\bf Case n=1}
  Using (\ref{F1(u|x)}) we get (\ref{P mu>0 h<mu}) since,  for $x\in[0,h)$, we have 
\begin{equation}
    P(1|x)=1-\frac{1}{2}e^{\frac{h-\mu-x}{\sigma}}
\end{equation}
from which we recognize $\eta_1$ and $\beta_1$.

{\bf Case n=2}
Conditioning on the position reached at time $1$ and observing that $h \leq \mu \leq \mu+x $ we get

\begin{align}
 P(2|x)  &=\int_0^h\left( \frac{e^{\frac{y-\mu-x}{\sigma}}}{2\sigma} + \frac{e^{-\frac{\mu+x}{\sigma}}}{2}\delta(y)\right)\left( 1-\frac{1}{2}e^{\frac{h-\mu-y}{\sigma}}\right)dy \nonumber\\ 
 &=\left( \frac{e^{\frac{h-\mu}{\sigma}}}{2}-\left(\frac{h}{4\sigma}+\frac{1}{4}\right)e^{\frac{h-2\mu}{\sigma}}    \right) e^{-\frac{x}{\sigma}}
\end{align}

where we recognize the coefficient $\beta_2$ in (\ref{beta2 Teor4.2}).

{\bf Case $n>2$}

Conditioning on the position reached at time $1$ we get
\begin{subequations}
\begin{align}
 P(n+1|x)=&\int_0^h\left( \frac{e^{\frac{y-\mu-x}{\sigma}}}{2\sigma} + \frac{e^{\frac{-\mu-x}{\sigma}}}{2}\delta(y)\right)\beta_{n}  e^{\frac{-y}{\sigma}}dy \\ 
   =&\left[\beta_{n}\left(\frac{1}{2}+\frac{h}{2\sigma}\right)e^{-\frac{\mu}{\sigma}} \right]e^{-\frac{x}{\sigma}}
\end{align}
\end{subequations}
The coefficients become
\begin{subequations}
\begin{align}
      \beta_{n+1}=&\beta_{n}\left(\frac{1}{2}+\frac{h}{2\sigma}\right)e^{-\frac{\mu}{\sigma}} \\
      =&\beta_2\left(\frac{1}{2}+\frac{h}{2\sigma}\right)^{n-1}e^{-\frac{(n-1)\mu}{\sigma}}\\
      =&\frac{e^{\frac{h-n\mu}{\sigma}}}{2}\left(\frac{1}{2}+\frac{h}{2\sigma}\right)^{n-1}-\frac{e^{\frac{h-(n+1)\mu}{\sigma}}}{2}\left(\frac{1}{2}+\frac{h}{2\sigma}\right)^{n}
\end{align}
\end{subequations}

\end{proof}

\subsection{Proof of Theorem \ref{Theorem FPT mu<0 }}
\begin{proof}
  {\bf Case n=1}\\
  Mimicking the proof for $n=1$ in Theorem \ref{Theorem FPT mu>0 h>mu}, with $0\leq x <h< h-\mu$, we have 
\begin{align}
	P(1|x)=\frac{1}{2}e^{\frac{\mu-h+x}{\sigma}}
\end{align}
where we recognize $\alpha_1^{(1,0)}$ as (\ref{alpha_1 T3}).

{\bf Case n}\\
 We assume that (\ref{P_n mu<0}) holds for $n$, and we show that it holds for $n+1$, for $n\geq 1$. 

Conditioning on the position reached at time $1$, using (\ref{smoluk}) we have
\begin{equation}
    P(n+1|x)=\sum_{i=1}^{h_n}\int_0^hP_i(n|y)f_Z(y,1|x,0)dy
\end{equation}

The structure of $P_i(n+1|x)$ changes according to $i=1$  or $1<i\leq h_{n+1}$.

If $i=1$ we have $x<-\mu$ and  
    \begin{align}
     P_1(n+1|x) &=\left\{\sum_{j=1}^{h_n}\int_{I_n^j}\frac{e^{-\frac{y-x-\mu}{\sigma}}}{2\sigma}P_j(n|y)dy+\left(1-\frac{e^{\frac{\mu+x}{\sigma}}}{2}\right)P_1(n|0)\right\}\mathbbm{1}_{I_{n+1}^1}(x)\nonumber\\
     &=\left\{\frac{e^{\frac{x+\mu}{\sigma}}}{{2^{n+1}\sigma^n}}\sum_{j=1}^{h_n}\left\{ \sum_{k=0}^{n-1}\int_{I_n^j} \left[\alpha_{n}^{(j,k)}(y+(n-1)\mu)^k +\int_{I_n^j}\beta_{n}^{(j,k)}(y+(n-1)\mu)^k e^{-\frac{2y}{\sigma}}\right]+2^n \sigma^{n-1}\eta_{n}^{j}e^{-\frac{y}{\sigma}}\right\}dy\right.\nonumber\\
     &\left.+\left(1-\frac{e^{\frac{\mu+x}{\sigma}}}{2}\right)P_1(n|0)\right\}\mathbbm{1}_{I_{n+1}^1}(x)\nonumber\\
     &=\left\{\frac{e^{\frac{x+\mu}{\sigma}}}{{2^{n+1}\sigma^n}}\left\{ \sum_{j=1}^{h_n}K_n^j -K_n^0 \right\} + \frac{K_n^0}{2^n\sigma^n}\right\}\mathbbm{1}_{I_{n+1}^1}(x)
     \end{align}
and we recognize the parameters (\ref{coeff N4.4 i=1}).

If $i>1$ we have $x \geq -\mu$ and 
\begin{align}\label{P_n+1^i_prima}
  	P_i(n+1|x)&= \left\{
  	P_1(n|0) \frac{e^{-\frac{\mu+x}{\sigma}}}{2}+\sum_{j=1}^{i^*-1}\int_{I^j_n} P_j(n|y)\frac{e^{\frac{y-\mu-x}{\sigma}}}{2\sigma}dy
	+\int_{-(i^*-1)\mu}^{\mu+x}  P_{i^*}(n|y)\frac{e^{\frac{y-\mu-x}{\sigma}}}{2\sigma}dy \right.\\	
	&+\left. \int_{\mu+x}^{-i^*\mu} P_{i^*}(n|y)\frac{e^{\frac{-y+\mu+x}{\sigma}}}{2\sigma}dy
	+\sum_{j=i^*+1}^{h_n}\int_{I^j_n} P_j(n|y)\frac{e^{\frac{-y+\mu+x}{\sigma}}}{2\sigma}dy  \right \} \mathbbm{1}_{I^{i}_{n+1}}(x). \nonumber
\end{align}
Now, using  the inductive hypothesis (\ref{P_ni mu<0}), we have

	\begin{align}
		P_i(n+1|x)&= \left\{\frac{e^{-\frac{\mu+x}{\sigma}}}{2}
		\left[\sum_{k=0}^{n-1} \frac{1}{2^n \sigma^{n-1}}\left[\alpha_{n}^{(1,k)}((n-1)\mu)^k +\beta_{n}^{(1,k)}((n-1)\mu)^k\right] + \eta_n^1\right]
		\right.\\
		&+\sum _{j=1}^{i^*-1} \frac{1}{2^n \sigma^{n-1}} \int_{I_n^{j}} \frac{e^{\frac{(y-\mu-x)}{\sigma}}}{2\sigma}\left\{ \sum_{k=0}^{n-1} \left[\alpha_{n}^{(j,k)}(y+(n-1)\mu)^k e^{\frac{y}{\sigma}}+\beta_{n}^{(j,k)}(y+(n-1)\mu)^k e^{-\frac{y}{\sigma}}\right]+\eta_{n}^{j}2^n \sigma^{n-1}\right\} dy\nonumber\\
		&+\frac{1}{2^n \sigma^{n-1}} \int_{-(i^*-1)\mu}^{x+\mu} \frac{e^{\frac{(y-\mu-x)}{\sigma}}}{2\sigma} \left\{\sum_{k=0}^{n-1} \left[\alpha_{n}^{(i^*,k)}(y+(n-1)\mu)^k e^{\frac{y}{\sigma}}+\beta_{n}^{(i^*,k)}(y+(n-1)\mu)^k e^{-\frac{y}{\sigma}}\right] +\eta_{n}^{i^*}2^n \sigma^{n-1}\right\}dy\nonumber\\
		&+\frac{1}{2^n \sigma^{n-1}} \int_{x+\mu}^{
  -i^*\mu} \frac{e^{\frac{-(y-\mu-x)}{\sigma}}}{2\sigma} \left\{\sum_{k=0}^{n-1}\left[\alpha_{n}^{(i^*,k)}(y+(n-1)\mu)^k e^{\frac{y}{\sigma}}+\beta_{n}^{(i^*,k)}(y+(n-1)\mu)^k e^{-\frac{y}{\sigma}}\right]+\eta_{n}^{i^*}2^n \sigma^{n-1}\right\} dy\nonumber\\
		&\left.+\sum _{j=i^*+1}^{h_n} \frac{1}{2^n \sigma^{n-1}} \int_{I_n^j} \frac{e^{\frac{-(y-\mu-x)}{\sigma}}}{2\sigma}\left\{\sum_{k=0}^{n-1} \left[\alpha_{n}^{(j,k)}(y+(n-1)\mu)^k e^{\frac{y}{\sigma}}+\beta_{n}^{(j,k)}(y+(n-1)\mu)^k e^{-\frac{y}{\sigma}}\right]+\eta_{n}^{j}2^n \sigma^{n-1}\right\} dy \right \} \nonumber\\&\mathbbm{1}_{I_{n+1}^i}(x)\nonumber
	\end{align}

Computing the integrals, we get

	\begin{align}\label{espr gamma}
	P_i(n+1|x)&=\left\{\frac{e^{-\frac{\mu+x}{\sigma}}}{2}\frac{K_n^0}{(2\sigma)^n}+
	\frac{e^{-\frac{(x+\mu)}{\sigma}}}{2^{n+1} \sigma^n} \sum _{j=1}^{i^*-1} K_n^j\right.\nonumber\\
	&+\frac{e^{-\frac{(x+\mu)}{\sigma}}}{2^{n+1} \sigma^n}\left\{\sum_{k=0}^{n-1} \left[\alpha_{n}^{(i^*,k)}\frac{e^{-\frac{2(n-1)\mu}{\sigma}}\sigma^{k+1} (-1)^k}{2^{k+1}}\left.\Gamma\left(1+k,\frac{-2(y+(n-1)\mu)}{\sigma}\right)\right|_{-(i^*-1)\mu}^{x+\mu} \right.\right.\nonumber\\
	&\left.+\beta_{n}^{(i^*,k)}\left.\frac{(y+(n-1)\mu)^{k+1}}{k+1}\right|_{-(i^*-1)\mu}^{x+\mu}\right]
	\left.\left. +\eta_{n}^{i^*}(2\sigma)^n e^{\frac{y}{\sigma}}\right|_{-(i^*-1)\mu}^{x+\mu} \right\} \nonumber\\
	&+\frac{e^{\frac{(x+\mu)}{\sigma}}}{2^{n+1} \sigma^n} \left\{\sum_{k=0}^{n-1} \left[\alpha_{n}^{(i^*,k)}\left.\frac{(y+(n-1)\mu)^{k+1}}{k+1}\right|_{\mu+x}^{-i^*\mu}\right.\right.\nonumber\\
	&-\beta_{n}^{(i^*,k)}\frac{e^{\frac{2(n-1)\mu}{\sigma}}\sigma^{k+1} }{2^{k+1}}\left.\left.\Gamma\left(1+k,\frac{2(y+(n-1)\mu)}{\sigma}\right)\right|_{x+\mu}^{-i^*\mu} \right]
	\left.\left. -\eta_{n}^{i^*}(2\sigma)^n e^{-\frac{y}{\sigma}}\right|^{-i^*\mu}_{x+\mu} \right\} \nonumber\\ 
	&\left.+\frac{e^{\frac{(x+\mu)}{\sigma}}}{2^{n+1} \sigma^{n}} \sum _{j=i^*+1}^{h_n} K_n^j \right\} \mathbbm{1}_{I_{n+1}^i}(x)
	\end{align}

Expanding the incomplete Gamma function (\ref{gammaExpansion}) we get
\begin{align}\label{espr sviluppo}
	P_i(n+1,x)&=\left\{\frac{e^{-\frac{\mu+x}{\sigma}}}{2}\frac{K_n^0}{(2\sigma)^n}+
	\frac{e^{-\frac{(x+\mu)}{\sigma}}}{2^{n+1} \sigma^n} \sum _{j=1}^{i^*-1} K_n^j\right.\nonumber\\
	&+\frac{e^{-\frac{(x+\mu)}{\sigma}}}{2^{n+1} \sigma^n}\left\{\sum_{k=0}^{n-1} \left[\alpha_{n}^{(i^*,k)}\frac{\sigma^{k+1} (-1)^k}{2^{k+1}}\left.k! e^{\frac{2y}{\sigma}} \sum_{r=0}^k \frac{(-2)^r}{r! \sigma^r}(y+(n-1)\mu)^r\right|_{-(i^*-1)\mu}^{x+\mu}\right.\right.\nonumber\\
	&+\left.\beta_{n}^{(i^*,k)}\left.\frac{(y+(n-1)\mu)^{k+1}}{k+1}\right|_{-(i^*-1)\mu}^{x+\mu}
	\right] \left.\left.+\eta_{n}^{i^*}(2\sigma)^n e^{\frac{y}{\sigma}}\right|_{-(i^*-1)\mu}^{x+\mu} \right\} \nonumber\\
	&+\frac{e^{\frac{(x+\mu)}{\sigma}}}{2^{n+1} \sigma^n}\left\{\sum_{k=0}^{n-1} \left[\alpha_{n}^{(i^*,k)}\left.\frac{(y+(n-1)\mu)^{k+1}}{k+1}\right|_{\mu+x}^{-i^*\mu}
	-\beta_{n}^{(i^*,k)}\frac{\sigma^{k+1} }{2^{k+1}}\left.k! e^{\frac{-2y}{\sigma}} \sum_{r=0}^k \frac{(2)^r}{r! \sigma^r}(y+(n-1)\mu)^r\right|^{-i^*\mu}_{x+\mu}\right]\right.\nonumber\\ 
	 &\left.\left. -\eta_{n}^{i^*}(2\sigma)^n e^{-\frac{y}{\sigma}}\right|^{-i^*\mu}_{x+\mu} \right\} \left.+\frac{e^{\frac{(x+\mu)}{\sigma}}}{2^{n+1} \sigma^n} \sum _{j=i^*+1}^{\ell_n} K_n^j  \right\} \mathbbm{1}_{I_{n+1}^i}(x)
\end{align}
Using (\ref{espr gamma}) and (\ref{espr sviluppo}) we recognize the coefficients (\ref{coeff N4.4}).
\end{proof}

\subsection{Proof of Corollary \ref{Theorem FPT mu<0 0<h< -mu}}
\begin{proof}
We could show that (\ref{P_n mu<0}) in Theorem \ref{Theorem FPT mu<0 } simplifies to  (\ref{P_n mu<0 h<-mu}). However, it is easier to proceed by induction. 
 
 {\bf Case n=1}
  Since $0<x\leq h \leq -\mu$ implies $x+\mu\leq 0$, using (\ref{F1(u|x) mu<-x}), we have
\begin{align}\label{stop1b}
	P(1|x)=	\frac{1}{2}e^{\frac{\mu-h+x}{\sigma}}
\end{align}

{\bf Case n}

Using the inductive hypothesis as in the previous theorems, conditioning on the position reached at time $1$ we get

\begin{align}
 P(n+1|x) &=\int_0^h\left[ \frac{e^{\frac{-y+\mu+x}{\sigma}}}{2\sigma} + \left(1-\frac{e^{\frac{\mu+x}{\sigma}}}{2}\right)\delta(y)\right]\left(\alpha_ne^{\frac{y}{\sigma}}+\eta_n\right)dy \nonumber\\ 
 &=h\frac{e^{\frac{\mu+x}{\sigma}}}{2\sigma}\alpha_n +\alpha_n-\frac{e^{\frac{\mu+x}{\sigma}}}{2}\alpha_n +\frac{e^{\frac{\mu+x}{\sigma}}}{2}\left(1-e^{-\frac{h}{\sigma}}\right)\eta_n+ \left(1-\frac{e^{\frac{\mu+x}{\sigma}}}{2}\right)\eta_n
\end{align}
Collecting the terms in $e^{\frac{x}{\sigma}}$  we find the recursive relations (\ref{alpha T4}) and (\ref{eta T4}).
\end{proof}

\subsection{Proof of Theorem \ref{Theorem FPT mu=0}}
\begin{proof}
We proceed by induction.

{\bf Case n=1}
 Since $0<x\leq h$, using (\ref{F1(u|x)}) we have
\begin{equation}
    P(1|x)=\frac{1}{2}e^{\frac{x-h}{\sigma}}=\frac{1}{2}\alpha_1^{(1,0)}e^{\frac{x}{\sigma}}
\end{equation}

{\bf Case n}
Using the inductive hypothesis as in the previous theorems, conditioning on the position reached at time 1 we get
\begin{align}
P(n+1|x)&=\int_0^h \frac{1}{2^n \sigma^{n-1}}\left[\sum_{j=0}^{n-1} \alpha_{n}^{(1,j)}y^j e^{\frac{y}{\sigma}}+\sum_{j=0}^{n-2} \beta_{n}^{(1,j)}y^j e^{-\frac{y}{\sigma}}\right]\frac{e^{-\frac{|y-x|}{\sigma}}}{2\sigma} dy+P(n|0)\frac{e^{-\frac{x}{\sigma}}}{2}\nonumber\\
&=\frac{e^{\frac{-x}{\sigma}}}{2^{n+1}\sigma^{n}}\left[\sum_{j=0}^{n-1} \alpha_{n}^{(1,j)}\frac{\sigma^{j+1} (-1)^j}{2^{j+1}}\left.\Gamma\left(1+j,\frac{-2y}{\sigma}\right)\right|_0^x+\sum_{j=0}^{n-2}\beta_{n}^{(1,j)}\left.\frac{y^{j+1}}{j+1}\right|_0^x\right]\nonumber\\
&+\frac{e^{\frac{x}{\sigma}}}{2^{n+1}\sigma^{n}}\left[\sum_{j=0}^{n-1} \alpha_{n}^{(1,j)}\left.\frac{y^{j+1}}{j+1}\right|_x^h+\sum_{j=0}^{n-2}\beta_{n}^{(1,j)}\frac{\sigma^{j+1}}{2^{j+1}}\left.\Gamma\left(1+j,\frac{2y}{\sigma}\right)\right|_x^h\right]+P(n|0)\frac{e^{-\frac{x}{\sigma}}}{2}
\end{align}
Expanding the incomplete Gamma function (\ref{gammaExpansion}) we recognize the coefficients (\ref{alfa beta mu=0})

\end{proof}

\section{Sensitivity analysis}
Here we apply the theorems presented in  Section \ref{Sect. distribution} to investigate the shapes of the distribution of $W_n$ for finite times, emphasizing the variety of the behaviours as parameters change.
\begin{figure}
\centering
\includegraphics[height= 8cm, keepaspectratio]{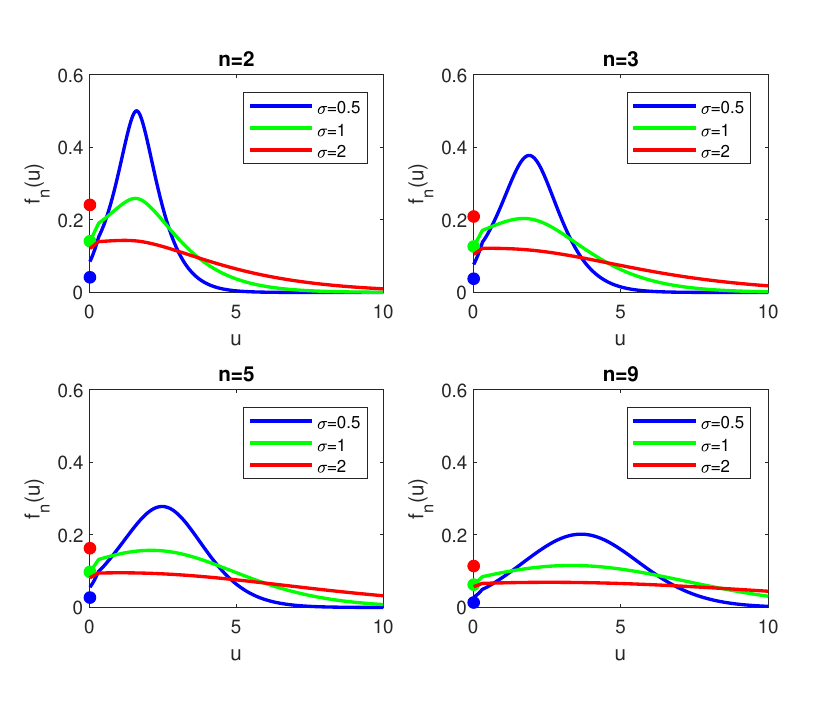}
\caption{Density functions of the Lindley process with $\mu=0.3$ and starting position $x=1$ for different values of $n=[2,3,5,9]$ and $\sigma$ ($\sigma=0.5$ blue, $\sigma=1$ green and $\sigma=2$ red). All these graphs have been obtained using Theorem \ref{Theorem Position mu>=0}.}\label{Fig posizione sigma mu0.3}
\end{figure}

\begin{figure}
\centering
\includegraphics[height= 8cm, keepaspectratio]{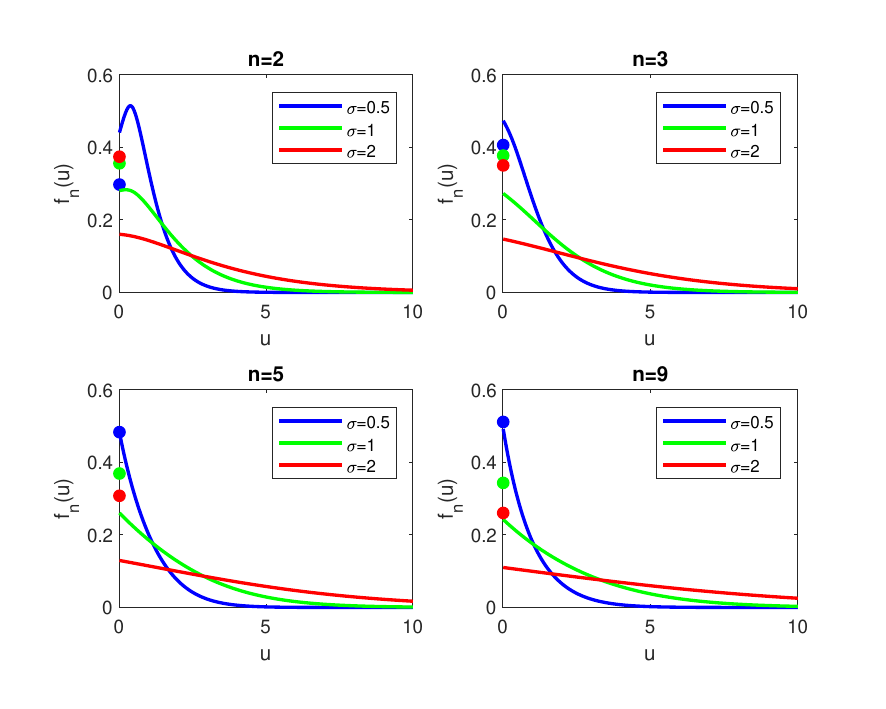}
\caption{Density functions of the Lindley process with $\mu=-0.3$ and starting position $x=1$ for different values of $n=[2,3,5,9]$ and $\sigma$ ($\sigma=0.5$ blue, $\sigma=1$ green and $\sigma=2$ red). All these graphs have been obtained using Theorem \ref{Theorem Position -x<mu<0}.}\label{Fig posizione sigma mu-0.3}
\end{figure}

In Figure \ref{Fig posizione sigma mu0.3} the density functions of the Lindley process with $\mu=0.3$ and starting position $x=1$ for different values of $n$ and $\sigma$ is shown. We can see that, as $n$ increases, the density flattens out, the variance increases
and the maximum of the density moves toward higher values. 
As $\sigma$ increases, the density flattens out but keeping fixed the position of the maximum.
The discrete part of the distribution, represented by a coloured dot on the $y$-axis decreases as $n$ increases. Note also that only when $\sigma=1$ we observe the continuity between the continuous
and discrete part of the density, this is due to Corollary \ref{Cor c_n+1}.

A different behaviour appears when the shift term $\mu$ is negative. In 
Figure \ref{Fig posizione sigma mu-0.3} the density functions of the Lindley process with the same parameters of Figure \ref{Fig posizione sigma mu0.3} and $\mu=-0.3$ is shown. We can see that, as $n$ increases, the density converges to a stationary distribution as expected from the theory. The interesting remark is that such convergence appears for reasonable small values of $n$. Furthermore, as $\sigma$ increases, the density flattens out and the variance of $W_n$ increases, as expected since we increased the variability of the process.

\begin{figure}
\centering
\includegraphics[height= 8cm, keepaspectratio]{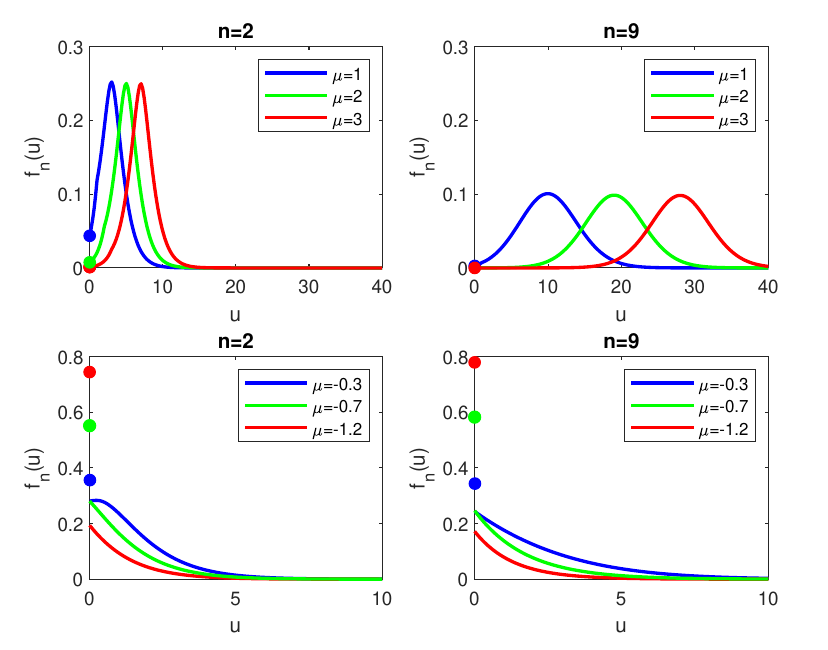}
\caption{Density functions of the Lindley process with $\sigma=1$ and starting position $x=1$ for different values of $n=[2,9]$ and $\mu$ (first line: $\mu=1$ blue, $\mu=2$ green and $\mu=3$ red; second line: $\mu=-0.3$ blue, $\mu=-0.7$ green and $\mu=-1.2$ red). The graphs of the first line  have been obtained using Theorem \ref{Theorem Position mu>=0}. In the second line for the graphs corresponding to $\mu=-1.2$ we used Theorem \ref{Theorem Position mu<=-x} while for the other cases we used Theorem \ref{Theorem Position -x<mu<0}.}\label{Fig posizione mu sigma1}
\end{figure}

Figure \ref{Fig posizione mu sigma1} shows the density function of the Lindley process with $\sigma=1$, starting position $x=1$ for different values of $n$ and $\mu$. We notice that, if $\mu$ is positive, the increase of $\mu$ determines very similar shape of the density but shifted while, if $\mu$ is negative, as $\mu$ decreases the density concentrates more and more in zero.
Another interesting feature concerns the mass in $u=0$. For positive $\mu$, it decreases as $\mu$ increases, indeed the trajectories move quickly away from zero. While, for negative $\mu$ we observe the opposite behavior due to the fact that the trajectories are pushed towards zero and it becomes increasingly difficult to leave zero. 

\begin{figure}
\centering
\includegraphics[height= 7cm, keepaspectratio]{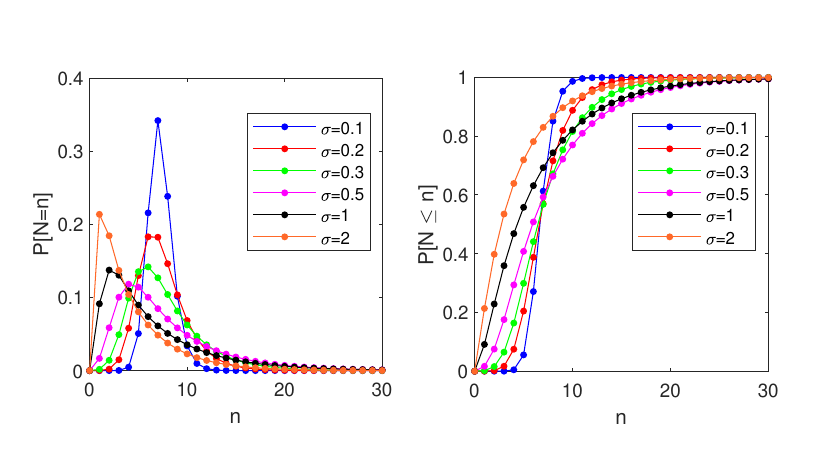}
\caption{Distribution and cumulative distribution function of the FET of the Lindley process originated in $x=1$ with $\mu=0.3$ and for different values of $\sigma$ ($\sigma=0.1$ blue, $\sigma=0.2$ red, $\sigma=0.3$ green, $\sigma=0.5$ magenta, $\sigma=1$ black and $\sigma=2$ orange). All these graphs have been obtained using Theorem \ref{Theorem FPT mu>0 h>mu}. In the plots we connected the probabilities to facilitate reading.}\label{Fig FETsigma mu0.3}
\end{figure}

As far as the FET $N_x$ are concerned, in Figure \ref{Fig FETsigma mu0.3} we illustrate the behaviour of the probability distribution function $P(n|x)$ and of its cumulative with starting position $x=1$, boundary $h=3$, $\mu=0.3$ and for different values of $\sigma$. 
We see a different behaviour for small or large values of the dispersion parameter $\sigma$. Indeed, when $\sigma \leq 0.5$, $P(n|x=1)$ has a maximum whose abscissa decreases as $\sigma$ decreases while it increases when $\sigma>0.5$. This behaviour has an immediate interpretation by observing that almost deterministic crossings determine an high peak when $\sigma$ is small; as far as $\sigma$ increases, the variability of the increments facilitates the crossing that also happens for smaller times and the probabilistic mass starts to increase sooner. 

\begin{figure}
\centering
\includegraphics[height= 9cm, keepaspectratio]{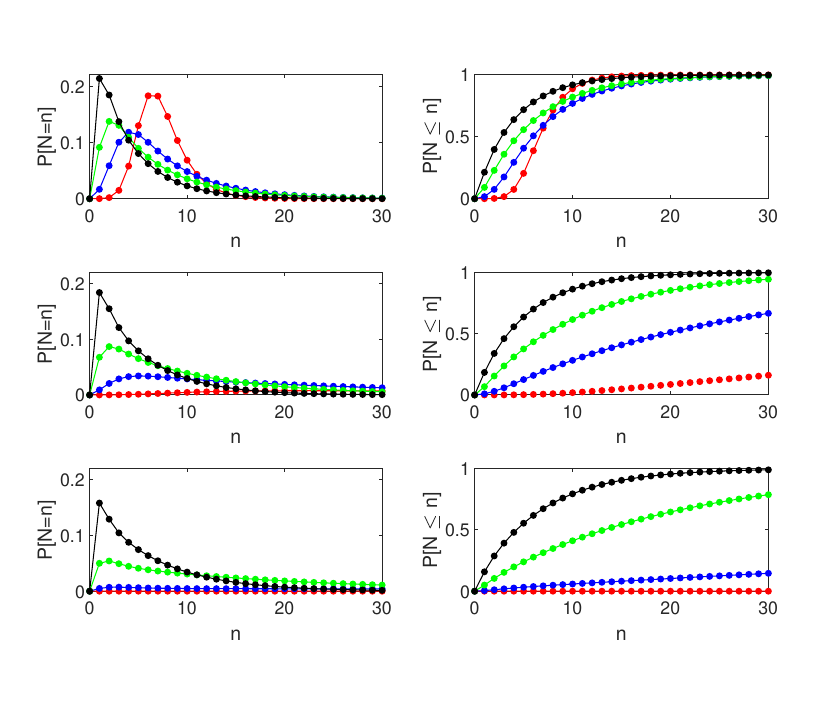}
\caption{Distribution and cumulative distribution function of the FET of the Lindley process originated in $x=1$ with $\mu=0.3$ (first row), $\mu=0$ (second row) and $\mu=-0.3$ (third row). For each choice of $\mu$ a comparison between different values of $\sigma$ is performed ($\sigma=0.2$ red, $\sigma=0.5$ blue, $\sigma=1$ green, $\sigma=2$ black). All these graphs have been obtained using Theorem \ref{Theorem FPT mu>0 h>mu} (first row), Theorem \ref{Theorem FPT mu=0} (second row) and Theorem \ref{Theorem FPT mu<0 } (third row). The probabilities in the plots are connected to facilitate reading.}\label{Fig FETsigma mu}
\end{figure}

In Figure \ref{Fig FETsigma mu} we investigate how this behaviour evolves when $\mu$ decreases. In particular we compare the shapes for different values of $\mu=[-0.3, 0, 0.3]$. Observe that the abscissa of the maximum of the distribution decreases as $\sigma$ increases. Concerning the corresponding ordinate we observe different behaviours depending on the sign of the parameter $\mu$. Indeed, for positive $\mu$ we have the features already noted in Figure \ref{Fig FETsigma mu0.3}. This fratures are no more observed when $\mu\leq 0$ because here the deterministic crossing is no more possible and crossings are determined only by the noise.  

\section*{Acknowledgement}
We are grateful to Professor I. Meilijson for his  useful suggestions and comments.

\end{document}